\documentclass[11pt]{article}
\usepackage{latexsym}
\usepackage[all]{xy}
\usepackage{amsmath}
\usepackage{amssymb}
\usepackage{amsfonts}
\usepackage[applemac]{inputenc}

\newtheorem{thm}{Theorem}[section]
\newtheorem{prop}[thm]{Proposition}
\newtheorem{lem}[thm]{Lemma}

\newtheorem{df}[thm]{Definition}
\newtheorem{cor}[thm]{Corollary}

\newtheorem{rmk}[thm]{Remark}

\newcommand  {\alg}   {\mathbf{Alg}_{\mathbb{C}}}
\newcommand  {\salg}   {\mathbf{SAlg}_{\mathbb{C}}}

\newcommand  {\ho}   {\mathrm{Ho}}
\newcommand  {\st}   {\mathbf{St}_{\mathbb{C}}}
\newcommand  {\dst}   {\mathbf{dSt}_{\mathbb{C}}}
\newcommand  {\ssets}   {\mathbf{SSets}}
\newcommand  {\tr}     {\mathrm{t}_{0}}

\begin{document}

\title{\textbf{Derived algebraic geometry, determinants of perfect complexes, and applications to obstruction theories for maps and complexes \bigskip \bigskip}}
\bigskip
\bigskip

\author{\bf{Timo Sch\"urg}\\ \small{Max-Planck-Institut f\"ur Mathematik} \\ \small{Bonn - Germany}
\and \bf{Bertrand To\"en}\\ \small{I3M UMR 5149} \\  \small{Universit\'e de Montpellier2 - France} \\ \small{Montpellier - France} \and \bf{Gabriele Vezzosi}\\ \small{Dipartimento di Sistemi ed Informatica} \\ \small{Sezione di Matematica} \\ \small{Universit\`a di Firenze}\\ \small{Firenze - Italy} \bigskip }

\date{June 2011}

\maketitle
\begin{abstract} We show how a quasi-smooth derived enhancement of a Deligne-Mumford stack $\mathcal{X}$ naturally endows $\mathcal{X}$ with a functorial perfect obstruction theory in the sense of Behrend-Fantechi. This result is then applied to moduli of maps and perfect complexes on a smooth complex projective variety.\\
For \emph{moduli of maps}, for $X=S$ an algebraic $K3$-surface, $g\in \mathbb{N}$, and $\beta \neq 0$ in $H_{2}(S,\mathbb{Z})$ a curve class, we construct a derived stack $\mathbb{R}\overline{\mathbf{M}}^{\textrm{red}}_{g,n}(S;\beta)$ whose truncation is the usual stack $\overline{\mathbf{M}}_{g,n}(S;\beta)$ of pointed stable maps from curves of genus $g$ to $S$ hitting the class $\beta$, and such that the inclusion $\overline{\mathbf{M}}_{g}(S;\beta)\hookrightarrow \mathbb{R}\overline{\mathbf{M}}^{\textrm{red}}_{g}(S;\beta)$ induces on $\overline{\mathbf{M}}_{g}(S;\beta)$ a perfect obstruction theory whose tangent and obstruction spaces coincide with the corresponding \emph{reduced} spaces of Okounkov-Maulik-Pandharipande-Thomas \cite{op, mp, mpt}. The approach we present here uses derived algebraic geometry and yields not only a full rigorous proof of the existence of a reduced obstruction theory - not relying on any result on semiregularity maps - but also a new global geometric interpretation. \\ 
We give two further applications to \emph{moduli of complexes}. For a $K3$-surface $S$ we show that the stack of simple perfect complexes on $S$ is smooth. This result was proved with different methods by Inaba (\cite{ina}) for the corresponding coarse moduli space. Finally, we construct a map from the derived stack of stable embeddings of curves (into a smooth complex projective variety $X$) to the derived stack of simple perfect complexes on $X$ with vanishing negative Ext's, and show how this map induces a morphism of the corresponding obstruction theories when $X$ is a Calabi-Yau threefold. \\
An important ingredient of our construction is a \emph{perfect determinant map} from the derived stack of perfect complexes to the derived stack of line bundles whose tangent morphism is, pointwise, Illusie's trace map for perfect complexes. We expect that this determinant map might be useful in other contexts as well.\end{abstract}


\bigskip

\tableofcontents

\section*{Introduction}

It is well known in Algebraic Geometry - e.g. in Gromov-Witten and Donaldson-Thomas theories - the importance of endowing a Deligne-Mumford moduli stack with a (perfect)  \emph{obstruction theory}, as defined in \cite{bf}: such an obstruction theory gives a \emph{virtual fundamental class} in the Chow group of the stack. If the stack in question is the stack of pointed stable maps to a fixed smooth projective variety (\cite{bm}), then integrating appropriate classes against this class produces all versions of Gromov-Witten invariants (\cite{be}).\\ Now, it is a distinguished feature of \emph{Derived Algebraic Geometry} (as exposed, e.g. in \cite{hagII}) that any quasi-smooth \emph{derived extension} of such a  stack $F$, i.e. a derived stack whose underived part or \emph{truncation} is the given stack $F$, induces a \emph{canonical} obstruction theory on $F$:  we have collected these results in \S \ref{derobs} below. A morphisms of derived stacks induces naturally a morphism between the induced obstruction theories - so that functoriality results like \cite[Prop. 5.10]{bf} or the so-called virtual pullback result in \cite{manolache} follow immediately. Moreover the functoriality of obstruction theories induced by morphisms of derived extensions is definitely richer than the usual one in \cite{bf}, that is restricted to special situations (e.g. \cite[Prop. 5.10]{bf}), and requires the axiomatics of compatible obstruction theories.  In other words, a suitable reformulation of a moduli problem in derived algebraic geometry, immediately gives us a canonical obstruction theory, in a completely geometric way, with no need of clever choices.  And, under suitable conditions, also the converse is expected to hold. \\ The present paper is an  application of this feature of derived algebraic geometry, \emph{not} relying on the mentioned conjectural general equivalence between a class of derived stacks and a class of underived stacks endowed with a properly structured obstruction theory. However, as a matter of fact, all the obstruction theories we are aware of indeed arise from derived extensions - and the cases covered in this paper simply add to this list. \\

In this paper we apply this ability of derived algebraic geometry in producing obstruction theories - functorial with respect to maps of derived stacks - to the cases of moduli of maps and moduli of perfect complexes on a complex smooth projective variety $X$.\\

\noindent \textbf{Moduli of maps.} For moduli of maps, we show how the standard obstruction theory yielding Gromov-Witten invariants comes from a natural derived extension of the stack of pointed stable maps to $X$. Then we concentrate on the first geometrically interesting occurrence of two different obstruction theories on a given stack, the stack $\overline{\mathbf{M}}_{g}(S;\beta)$ of stable maps of type $(g, \beta)$ to a smooth projective complex $K3$-surface $S$. The stack $\overline{\mathbf{M}}_{g}(S;\beta)$ has a \emph{standard} obstruction theory, yielding trivial Gromov-Witten invariants in the $n$-pointed case, and a so-called \emph{reduced} obstruction theory, first considered by Okounkov-Maulik-Pandharipande-Thomas (often abbreviated to O-M-P-T in the text), giving interesting - and extremely rich in structure - curve counting invariants in the $n$-pointed case (see \cite{pak3, mp, mpt}, and \S \ref{review} below, for a detailed review). In this paper we use derived algebraic geometry to give a construction of a global reduced obstruction theory on $\overline{\mathbf{M}}_{g}(S;\beta)$, and compare its deformation and obstruction spaces with those of Okounkov-Maulik-Pandharipande-Thomas. More precisely, we use a perfect determinant map form the derived stack of perfect complexes to the derived stack of line bundles, and exploit the peculiarities of the derived stack of line bundles on a $K3$-surface, to produce a derived extension $\mathbb{R}\overline{\mathbf{M}}^{\textrm{red}}_{g}(S;\beta)$ of $\overline{\mathbf{M}}_{g}(S;\beta)$. The derived stack $\mathbb{R}\overline{\mathbf{M}}^{\textrm{red}}_{g}(S;\beta)$ arises as the canonical homotopy fiber over the \emph{unique derived factor} of the derived stack of line bundles on $S$, so it is, in a very essential way, a purely derived geometrical object. We prove quasi-smoothness of $\mathbb{R}\overline{\mathbf{M}}^{\textrm{red}}_{g}(S;\beta)$, and this immediately gives us a global reduced obstruction theory on $\overline{\mathbf{M}}_{g}(S;\beta)$. Our proof is self contained (inside derived algebraic geometry), and does not rely on any previous results on semiregularity maps.\\

\noindent \textbf{Moduli of complexes.} We give two applications to moduli of perfect complexes on smooth projective varieties.
In the first one we show that the moduli space of simple perfect complexes on a $K3$-surface is smooth. Inaba gave a direct proof of this result in \cite{ina}, by generalizing methods of Mukai (\cite{muk}). Our proof is different and straightforward. We use the perfect determinant map, and the peculiar structure of the derived Picard stack of a $K3$-surface, to produce a derived stack of simple perfect complexes. Then we show that this derived stack is actually \emph{underived} (i.e trivial in the derived direction) and \emph{smooth}. The moduli space studied by Inaba is exactly the coarse moduli space of this stack.\\
In the second application, for $X$ an arbitrary smooth complex projective scheme $X$, we first construct a map $C$ from the derived stack $\mathbb{R}\overline{\mathbf{M}}_{g,n}(X)^{emb}$ consisting of pointed stable maps which are closed immersions, to the derived stack $\mathbb{R}\mathbf{Perf}(X)^{\textrm{si}, > 0}_{\mathcal{L}}$ of simple perfect complexes with no negative Ext's and fixed determinant $\mathcal{L}$ (for arbitrary $\mathcal{L}$). Then we show that, if $X$ is a \emph{Calabi-Yau threefold}, the derived stack $\mathbb{R}\mathbf{Perf}(X)^{\textrm{si}, > 0}_{\mathcal{L}}$ is actually \emph{quasi-smooth}, and use the map $C$ to compare (according to \S \ref{gwdt}) the canonical obstruction theories induced by the source and target derived stacks on their truncations. Finally, we relate this second applications to a baby, open version of the Gromov-Witten/Donaldson-Thomas conjectural comparison. In such a comparison, one meets two basic problems. The first, easier, one is in producing a map enabling one to compare the obstruction theories - and derived algebraic geometry, as we show in the open case, is perfectly suited for this (see \S \ref{derobs} and \S \ref{gwdt}). Such a comparison would induce a comparison (via a virtual pullback construction as in \cite[Thm 7.4]{sch}) between the corresponding \emph{virtual fundamental classes}, and thus a comparison between the GW and DT invariants. The second problem, certainly the most difficult one, is to deal with problems arising at the boundary of the compactifications. For this second problem, derived algebraic methods unfortunately do not provide at the moment any new tool or direction.\\

One of the main ingredients of all the applications given in this paper is the construction of a \emph{perfect determinant map}  $\det_{\textbf{Perf}}: \mathbf{Perf}\rightarrow \mathbf{Pic}$, where $\mathbf{Perf}$ is the stack of perfect complexes, $\mathbf{Pic}$ the stack of line bundles, and both are viewed as derived stacks (see \S 3.1 for details), whose definition requires the use of a bit of Waldhausen $K$-theory for simplicial commutative rings, and whose tangent map can be identified with Illusie's trace map of perfect complexes (\cite[Ch. 5]{ill}). We expect that this determinant map might be useful in other moduli contexts as well.\\

An important remark - especially for applications to Gromov-Witten theory - is that, in order to simplify the exposition, we have chosen to write the proofs only in the non-pointed case, since obviously no substantial differences except for notational ones are involved. The relevant statements are however given in both the unpointed and the $n$-pointed case.\\

To summarize, there are four main points in our paper.  The \emph{first} one is that derived algebraic geometry is a natural world where to get completely functorial obstruction theories. This is proved in \S 1, and explained through many examples in the rest of the text.\\ 
The \emph{second} main point concerns the application to reduced obstruction theory for stable maps to a $K3$ surface. More precisely, we give a rigorous proof of the existence of a global reduced obstruction theory on the stack of pointed stable maps to a $K3$-surface. The most complete, among the previous attempts in the literature, is the unpublished \cite[\S 2.2]{mp}, that only leads - after some further elaboration - to a uniquely defined  ``obstruction theory'' with target just the \emph{truncation in degrees} $\geq -1$ of the cotangent complex of the stack of maps (from a fixed domain curve). Moreover, in order to reach this, the authors invoke some results on semiregularity maps, whose validity does not seem to be completely satisfactorily established. Nevertheless, there is certainly a clean and complete description of the corresponding expected tangent and obstruction spaces in several papers (see \cite{op, mp, mpt}), and we prove that our obstruction theory has exactly such tangent and obstruction spaces. Our approach does not only establish rigorously such a reduced \emph{global} obstruction theory - with values in the \emph{full} cotangent complex of the stack of stable maps - but also endows such an obstruction theory with all the functoriality properties that sometimes are missing in the underived, purely obstruction-theoretic approach. This might prove useful in getting similar results in families or even relative to the whole moduli stack of $K3$-surfaces. Our definition of the reduced obstruction theory on the stack of stable maps to a $K3$-surface comes together with a clear (derived) geometrical picture - half of which is valid for \emph{any} smooth complex projective variety. This should be compared to the other existing, partial approaches, where the construction of obstruction spaces arises from local linear algebra manipulations whose global geometrical interpretation is a bit obscure. In order to get a satisfying global geometric picture underlying this reduced obstruction theory, we are forced to move to the world of derived algebraic geometry; so, in some sense, our picture describes and explains the geometry underlying those local computations. This is another example of what seems to be a fairly general principle: some constructions on moduli stacks, that happen to be ad hoc inside algebraic geometry, become canonical and gain a neater geometrical interpretation in derived algebraic geometry.\\
The \emph{third} main point of our paper is another application of our perfect determinant map. We show that the stack of simple perfect complexes on a $K3$-surface is smooth. In the \emph{fourth} and final application, for $X$ a Calabi-Yau threefold, we use the perfect determinant map and the functoriality of obstruction theories arising from derived extensions, to produce a map from a derived stack of stable maps to $X$ to the derived stack of simple perfect complexes with fixed determinant on $X$. We prove that the target derived stack is quasi-smooth, and show how this map induces, by functoriality, a comparison map between the associated obstruction theories.\\

Finally, let us observe that the emerging picture seems to suggest that most of the natural maps of complexes arising in moduli problems can be realized as tangent or cotangent maps associated to morphisms between appropriate derived moduli stacks. This suggestion is confirmed in the present paper for the standard obstruction theories associated to the stack of maps between a fixed algebraic scheme and a smooth projective target, to the stack of stable maps to a smooth projective scheme or to the Picard stack of a smooth projective scheme, for the trace map, the Atiyah class map, and the first Chern class map for perfect complexes (\cite[Ch. 5]{ill}), and for the map inducing O-M-P-T's reduced obstruction theory. \\

\noindent \textbf{Description of contents.} The first three sections and the beginning of the fifth are written for an \emph{arbitrary} smooth complex projective scheme $X$. We explain how a derived extension induces an obstruction theory on its truncation (\S 1), how to define the standard derived extensions of the Picard stack of $X$, and of the stack of stable maps to $X$ (\S 2), and finally define the perfect determinant map (\S 3). In section (\S 4), we \emph{specialize} to the case where $X=S$ is a smooth complex projective $K3$ surface. We first give a self-contained description of O-M-P-T's pointwise reduced tangent and obstruction spaces (\S 4.1). Then, by exploiting the features of the derived Picard stack of $S$ (\S 4.2), we define in \S 4.3 a derived extension $\mathbb{R}\overline{\mathbf{M}}^{\textrm{red}}_{g}(S;\beta)$ of the usual stack $\overline{\mathbf{M}}_{g}(S;\beta)$ of stable maps of type $(g,\beta\neq 0)$ to $S$, having the property that, for the canonical inclusion $j_{\textrm{red}}: \overline{\mathbf{M}}_{g}(S;\beta)\hookrightarrow \mathbb{R}\overline{\mathbf{M}}^{\textrm{red}}_{g}(S;\beta)$, the induced map $$j_{\textrm{red}}^{*}\mathbb{L}_{\mathbb{R}\overline{\mathbf{M}}^{\textrm{red}}_{g}(S;\beta)} \longrightarrow \mathbb{L}_{\overline{\mathbf{M}}_{g}(S;\beta)}$$ is a perfect obstruction theory with the same tangent and obstruction spaces as the reduced theory introduced by Maulik-Okounkov-Pandharipande-Thomas (\S 4.4, Theorem. \ref{quasismooth}). \\
In \S 5, for a complex smooth projective variety $X$, we define the derived stack $\mathbb{R}\overline{\mathbf{M}}_{g,n}(X)^{emb}$ of pointed stable maps to $X$ that are closed embeddings, the derived stack $\mathcal{M}_{X}\equiv \mathbb{R}\mathbf{Perf}(X)^{\textrm{si}, > 0}$ of simple perfect complexes on $X$ with vanishing negative Ext's, and the derived stack $\mathcal{M}_{X, \mathcal{L}}\equiv \mathbb{R}\mathbf{Perf}(X)^{\textrm{si}, > 0}_{\mathcal{L}}$ of simple perfect complexes on $X$ with vanishing negative Ext's and fixed determinant $\mathcal{L}$, and we define a morphism $C: \mathbb{R}\overline{\mathbf{M}}_{g,n}(X)^{emb} \longrightarrow \mathcal{M}_{X, \mathcal{L}}$. When $X$ is a $K3$-surface, we show that the truncation stack of $\mathcal{M}_{X}$ is smooth. When $X$ is a Calabi-Yau threefold, we prove that $\mathcal{M}_{X, \mathcal{L}}$ is quasi-smooth, and that the map $C$ induces a map between the obstruction theories on the underlying underived stacks.\\

In an Appendix we give a derived geometrical interpretation of the Atiyah class map and the first Chern class map for a perfect complex $\mathbf{E}$ on a scheme $Y$, by relating them to the tangent of the corresponding map $\varphi_{\mathbf{E}}: Y\rightarrow \mathbf{Perf}$; then we follow this reinterpretation to prove some properties used in the main text.\\

\bigskip

\noindent \textbf{Acknowledgments.} Our initial interest in the possible relationships between reduced obstruction theories and derived algebraic geometry was positively boosted by comments and questions by B. Fantechi, D. Huybrechts and R. Thomas. We are grateful to R. Pandharipande for pointing out a useful classical statement, and to H. Flenner for some important remarks. We especially thanks A. Vistoli for generously sharing his expertise on stable maps with us, and R. Thomas for his interest and further comments on this paper. \\
The second and third authors acknowledge financial support from the french ANR grant HODAG (ANR-09-BLAN-0151). \\

\bigskip

\noindent\textbf{Notations.} For background and basic notations in derived algebraic geometry we refer the reader to $\cite[Ch. 2.2]{hagII}$ and to the overview \cite[\S 4.2, 4.3]{seattle}. In particular, $\st$ (respectively, $\dst$) will denote the (homotopy) category of \emph{stacks} (respectively, of \emph{derived stacks}) on $\mathrm{Spec}\, \mathbb{C}$ with respect to the \'etale (resp., strongly \'etale) topology. We will most often omit the inclusion  functor $i: \st \rightarrow \dst $ from our notations, since it is fully faithful; its left adjoint, the truncation functor, will be denoted $\tr : \dst \rightarrow \st$ (\cite[Def. 2.2.4.3]{hagII}). In particular, we will write $\tr (F) \hookrightarrow F$ for the adjunction morphism $i\tr (F) \hookrightarrow F$. However recall that the inclusion functor $i$ does \emph{not} commute with taking internal HOM (derived) stacks nor with taking homotopy limits. All fibered products of derived stacks will be implicitly derived (i.e. they will be homotopy fibered products in the model category of derived stacks).\\ 
When useful, we will freely switch back and forth between (the model category of) simplicial commutative $k$-algebras and (the model category of) commutative differential non-positively graded $k$-algebras, where $k$ is a field of characteristic $0$ (\cite[App. A]{rhamloop}).\\ 
All complexes will be cochain complexes and, for such a complex $C^{\bullet}$, either $C_{\leq n}$ or $C^{\leq n}$ (depending on typographical convenience) will denote its good truncation in degrees $\leq n$. Analogously for either $C_{\geq n}$ or $C^{\geq n}$ (\cite[1.2.7]{we}).\\ To ease notation we will often  write $\otimes$ for the derived tensor product $\otimes^{\mathbb{L}}$, whenever no confusion is likely to arise.\\
$X$ will denote a smooth complex projective scheme while $S$ a smooth complex projective $K3$-surface.\\
As a purely terminological remark, for a given obstruction theory, we will call its \emph{deformation space} what is usually called its \emph{tangent space} (while we keep the terminology \emph{obstruction space}). We do this to avoid confusion with tangent spaces, tangent complexes or tangent cohomologies of related (derived) stacks.\\
We will often abbreviate the list of authors Okounkov-Maulik-Pandharipande-Thomas to O-M-P-T.\\

\section{Derived extensions, obstruction theories and their functoriality} \label{derobs}

We briefly recall here the basic observation that a derived extension of a given stack $\mathcal{X}$ induces an obstruction theory (in the sense of \cite{bf}) on $\mathcal{X}$, and deduce a richer functoriality with respect to the one known classically. Everything in this section is true over an arbitrary base ring, though it will be stated for the base field $\mathbb{C}$.\\

\subsection{Derived extensions induce obstruction theories} \label{derobss}

Let $\tr: \dst \rightarrow \st$ be the truncation functor between derived and underived stacks over $\mathbb{C}$ for the \'etale topologies (\cite[Def. 2.2.4.3]{hagII}). It has a left adjoint $i:\st \rightarrow \dst$ which is fully faithful (on the homotopy categories), and is therefore usually omitted from our notations. 

\begin{df} \label{derext}
Given a stack $\mathcal{X} \in \mathrm{Ho}(\st)$, a \emph{derived extension} of $\mathcal{X}$ is a derived stack $\mathcal{X}^{\mathrm{der}}$ together with an isomorphism $$\mathcal{X} \simeq \tr (\mathcal{X}^{\mathrm{der}}) .$$
\end{df}

\begin{prop}\label{inducedobstruction}
Let $\mathcal{X}^{\mathrm{der}}$ be a derived geometric stack which is a derived extension of the (geometric) stack $\mathcal{X}$. Then, the closed immersion $$j: \mathcal{X} \simeq \tr (\mathcal{X}^{\mathrm{der}}) \hookrightarrow \mathcal{X}^{\mathrm{der}}$$ induces a morphism $$j^{*} (\mathbb{L}_{\mathcal{X}^{\mathrm{der}}}) \longrightarrow \mathbb{L}_{\mathcal{X}}$$ which is 2-connective, i.e. its cone has vanishing cohomology in degrees $\geq -1$.
\end{prop}

\noindent \textbf{Proof.} The proof follows easily from the remark that if $A$ is a simplicial commutative $\mathbb{C}$-algebra and $A\rightarrow \pi_{0}(A)$ is the canonical surjection, then the cotangent complex $\mathbb{L}_{\pi_{0}(A)/A}$ is 2-connective, i.e. has vanishing cohomology in degrees $\geq -1$.
\hfill $\Box$ \\

The previous Proposition shows that a derived extension always induces an obstruction theory (whenever such a notion is defined by \cite[Def. 4.4]{bf}, e.g. when  $\mathcal{X}$ is a Deligne-Mumford stack). In particular, recalling that a derived stack is \emph{quasi-smooth} if its cotangent complex is perfect of amplitude in $[-1,0]$, we have the following result

\begin{cor} \label{cor} Let $\mathcal{X}^{\mathrm{der}}$ be a quasi-smooth derived Deligne-Mumford stack which is a derived extension of a (Deligne-Mumford)   stack $\mathcal{X}$. Then $$j^{*} (\mathbb{L}_{\mathcal{X}^{\mathrm{der}}}) \longrightarrow \mathbb{L}_{\mathcal{X}}$$ is a $[-1,0]$-perfect obstruction theory as defined in \cite[Def. 5.1]{bf}.
\end{cor}

\begin{rmk}\emph{We expect that also the converse is true, i.e. that given any stack $X$ locally of finite presentation over a field $k$, endowed with a map of co-dg-Lie algebroids $E \rightarrow \mathbb{L}_{X}$, there should exist a derived extension inducing the given obstruction theory. We will come back to this in a future work and will \emph{not} use it in the rest of this paper, although it should be clear that it was exactly such an expected result that first led us to think about the present work. Moreover, as a matter of fact, all the obstruction theories we are aware of indeed come from derived extensions, and the cases covered below simply add to this list.}
\end{rmk}
 
\subsection{Functoriality of deformation theories induced by derived extensions}\label{functobs}

If $f:X \rightarrow Y$ is a morphism of (Deligne-Mumford) stacks, and $o_X: E_X \rightarrow \mathbb{L}_{X}$ and $o_Y:E_Y\rightarrow \mathbb{L}_{Y}$ are ($[-1,0]$-perfect) obstruction theories, the classical theory of obstructions (\cite{bf}) does not provide in general a map 
$f^*E_Y \rightarrow E_X$ such that the following square $$\xymatrix{f^*E_Y \ar[r]^-{f^*(o_Y)} \ar[d] & f^*\mathbb{L}_Y \ar[d] \\ E_X \ar[r]^-{o_X} & \mathbb{L}_X}$$ is commutative (or commutative up to an isomorphism) in the derived category $\mathrm{D}(X)$ of complexes on $X$, where $ f^{*}\mathbb{L}_Y \rightarrow \mathbb{L}_X$ is the canonical map on cotangent complexes induced by $f$ (\cite[Ch. 2, (1.2.3.2)']{ill}). On the contrary, if $j_X: X\hookrightarrow \mathbb{R}X$ and $j_Y: Y \hookrightarrow \mathbb{R}Y$ are quasi-smooth derived (Deligne-Mumford) extensions of $X$ and $Y$, respectively, and $F:\mathbb{R}X \rightarrow \mathbb{R}Y$ is a morphism of derived stacks $$\xymatrix{X \ar[r]^-{\tr F} \ar@{^{(}->}[d] & Y \ar@{^{(}->}[d] \\ \mathbb{R}X \ar[r]^-{F} & \mathbb{R}Y},$$ then $j_X^*\mathbb{L}_{\mathbb{R}X} \rightarrow \mathbb{L}_{X}$ and $j_Y^*\mathbb{L}_{\mathbb{R}Y} \rightarrow \mathbb{L}_{Y}$ are ($[-1,0]$-perfect) obstruction theories by Cor. \ref{cor}, and \emph{moreover} there is indeed a canonical morphism of triangles in $\mathrm{D}(X)$ (we denote $\tr (F)$ by $f$) $$\xymatrix{f^*j_{Y}^*\mathbb{L}_{\mathbb{R}Y} \ar[r] \ar[d] & f^*\mathbb{L}_{Y} \ar[r] \ar[d] & f^*\mathbb{L}_{\mathbb{R}Y/Y} \ar[d] \\ j_{X}^*\mathbb{L}_{\mathbb{R}X} \ar[r] & \mathbb{L}_{X} \ar[r] & \mathbb{L}_{\mathbb{R}X/X} }$$ (see \cite[Prop. 1.2.1.6]{hagII} or \cite[Ch. 2, (2.1.1.5)]{ill}). This map relates the two induced obstruction theories and may be used to relate the corresponding virtual fundamental classes, too (when they exist, e.g. when $X$ and $Y$ are proper over $k$). We will not do this here since we will not need it for the results in this paper. However, the type of result we are referring to is the following

\begin{prop}\label{timo}\emph{\cite[Thm. 7.4]{sch}} Let $F:\mathbb{R}X \rightarrow \mathbb{R}Y$ be a quasi-smooth morphism between quasi-smooth Deligne-Mumford stacks, and $f:X\rightarrow Y$ the induced morphism on the truncations. Then, there is an induced virtual pullback (as defined in \emph{\cite{manolache}}) $f^!:A_*(Y)\rightarrow A_*(X)$, between the Chow groups of $Y$ and $X$, such that $f^!([Y]^\textrm{vir})= [X]^\textrm{vir}$, where $[X]^\textrm{vir}$ (respectively, $[Y]^\textrm{vir}$) is the virtual fundamental class (\emph{\cite{bf}}) on $X$ (resp., on $Y$) induced by the $[-1,0]$ perfect obstruction theory $j_X^*\mathbb{L}_{\mathbb{R}X} \rightarrow \mathbb{L}_{X}$ (resp., by  $j_Y^*\mathbb{L}_{\mathbb{R}Y} \rightarrow \mathbb{L}_{Y}$).
\end{prop}

\section{Derived stack of stable maps and derived Picard stack}

In this section we prove a correspondence between derived open substacks of a derived stack and open substacks of its truncation, and use it to construct the derived Picard stack $\mathbb{R}\mathbf{Pic}(X;\beta)$ of type $\beta \in H^{2}(X,\mathbb{Z})$, for any complex projective smooth variety $X$. After recalling the derived version of the stack of (pre-)stable maps to $X$, possibly pointed, the same correspondence will lead us to defining the derived stack $\mathbb{R}\overline{\mathbf{M}}_{g}(X;\beta)$ of stable maps of type $(g, \beta)$ to $X$ and its pointed version.\\

Throughout the section  $X$ will denote a smooth complex projective scheme, $g$ a nonnegative integer, $c_{1}$ a class in $H^{2}(X,\mathbb{Z})$ (which, for our purposes, may be supposed to belong to the image of $\mathrm{Pic}(X) \simeq H^{1}(X,\mathcal{O}^{*}_{X}) \rightarrow H^{2}(X,\mathbb{Z})$, i.e. belonging to $H^{1,1}(X) \cap H^{2}(X,\mathbb{Z})$), and $\beta \in H_{2}(X,\mathbb{Z})$ an effective curve class.\\

We will frequently use of the following

\begin{prop} \label{opens}Let $F$ be a derived stack and $\tr (F)$ its truncation. There is a bijective correspondence $$\phi_{F}: \{ \emph{Zariski open substacks of} \,\,\tr (F) \}  \longrightarrow \{ \emph{Zariski open derived substacks of} \,\, F \} .$$ For any Zariski open substack $U_0\hookrightarrow \tr (F)$, we have a homotopy cartesian diagram in $\dst$ $$\xymatrix{ U_{0} \ar@{^{(}->}[r] \ar[d] & \tr (F) \ar[d] \\ \phi_{F}(U_{0}) \ar@{^{(}->}[r] & F}$$ where the vertical maps are the canonical closed immersions.
\end{prop}

\noindent \textbf{Proof.} The statement is an immediate consequence of the fact that $F$ and $\tr (F)$ have the same topology (\cite[Cor. 2.2.2.9]{hagII}). More precisely, let us define $\phi _{F}$ as follows. If $U_{0}\hookrightarrow \tr (F)$ is an open substack, $\phi_{F}(U_{0})$ is the functor $$\salg \longrightarrow \ssets \, : \, A\longmapsto F(A)\times_{\tr (F)(\pi_{0} (A))} U_{0}(\pi_{0} (A))$$ where $F(A)$ maps to $\tr (F)(\pi_{0} (A))$ via the morphism (induced by the truncation functor $\tr$) $$F(A)\simeq \mathbb{R}\underline{\mathrm{Hom}}_{\dst}(\mathbb{R}\mathrm{Spec}(A), F)\longrightarrow \mathbb{R}\underline{\mathrm{Hom}}_{\st}(\tr (\mathbb{R}\mathrm{Spec}(A)), \tr (F))\simeq  \tr (F) (\pi_0 (A)).$$
The inverse to $\phi_{F}$ is simply induced by the truncation functor $\tr$.
\hfill $\Box$ \\

\subsection{The derived Picard stack}

\begin{df} The \emph{Picard stack} of $X/\mathbb{C}$ is the stack $$\mathbf{Pic}(X):= \mathbb{R}\mathrm{HOM}_{\st}(X,B\mathbb{G}_{m}).$$ The \emph{derived Picard stack} of $X/\mathbb{C}$ is the derived stack $$\mathbb{R}\mathbf{Pic}(X):= \mathbb{R}\mathrm{HOM}_{\dst}(X,B\mathbb{G}_{m}).$$
\end{df}

By definition we have a natural isomorphism $\tr (\mathbb{R}\mathbf{Pic}(X)) \simeq \mathbf{Pic}(X)$ in $\ho (\dst)$. Note that even though $\mathbf{Pic}(X)$ is smooth, it is \emph{not} true that $\mathbb{R}\mathbf{Pic}(X) \simeq \mathbf{Pic}(X)$, if $\mathrm{dim}(X)>1$; this can be seen on tangent spaces since $$\mathbb{T}_{L} \mathbb{R}\mathbf{Pic}(X) \simeq \textrm{C}^{\bullet}(X, \mathcal{O}_{X})[1]$$ for any global point $x_{L}:\mathrm{Spec}(\mathbb{C}) \rightarrow \mathbb{R}\mathbf{Pic}(X) $ corresponding to a line bundle $L$ over $X$. \\

Given $c_1 \in H^{2}(X,\mathbb{Z})$, we denote by $\mathbf{Pic}(X;c_1)$ the open substack of $\mathbf{Pic}(X)$ classifying line bundles with first Chern class $c_1$. More precisely, for any $R \in \alg$,  let us denote by $\textrm{Vect}_{1}(R;c_1)$ the groupoid of line bundles $L$ on $\mathrm{Spec}(R)\times X$ such that, for any point $x:\mathrm{Spec}(\mathbb{C}) \rightarrow \mathrm{Spec}(R)$ the pullback line bundle $L_{x}$ on $X$ has first Chern class equal to $c_1$. Then, $\mathbf{Pic}(X;c_1)$ is the stack: $$\alg \longrightarrow \mathbf{SSets} \, : \, R \longmapsto \mathrm{Nerve}(\textrm{Vect}_{1}(R;c_1))$$ where $\mathrm{Nerve}(C)$ is the nerve of the category $C$.\\
Note that we have $$\mathbf{Pic}(X)=\coprod_{c_1 \in H^{2}(X,\mathbb{Z})} \mathbf{Pic}(X;c_1). $$

\begin{df}  Let $c_1 \in H^{2}(X,\mathbb{Z})$. The \emph{derived Picard stack of type} $c_1$  of $X/\mathbb{C}$ is the derived stack $$\mathbb{R}\mathbf{Pic}(X;c_1):= \phi_{\mathbb{R}\mathbf{Pic}(X)}(\mathbf{Pic}(X;c_1)).$$
\end{df}

In particular, we have a natural isomorphism $\tr (\mathbb{R}\mathbf{Pic}(X;c_1)) \simeq \mathbf{Pic}(X;c_1)$, and a homotopy cartesian diagram in $\dst$ $$\xymatrix{ \mathbf{Pic}(X;c_1) \ar@{^{(}->}[r] \ar[d] & \mathbf{Pic}(X) \ar[d] \\ \mathbb{R}\mathbf{Pic}(X;c_1) \ar@{^{(}->}[r] & \mathbb{R}\mathbf{Pic}(X)}$$

\subsection{The derived stack of stable maps} \label{derivedstable}
We recall from \cite[4.3 (4.d)]{seattle} the construction of the derived stack $\mathbb{R}\mathbf{M}^{\textrm{pre}}_{g}(X)$ (respectively, $\mathbb{R}\mathbf{M}^{\textrm{pre}}_{g,n}(X)$) of prestable maps (resp., of $n$-pointed prestable maps) of genus $g$ to $X$, and of its open derived substack $\mathbb{R}\overline{\mathbf{M}}_{g}(X)$ (respectively, $\mathbb{R}\overline{\mathbf{M}}_{g,n}(X)$) of stable maps (resp., of $n$-pointed stable maps) of genus $g$ to $X$. Then we move to define the derived version of the stack of (pointed) stable maps of type $(g,\beta)$ to $X$. \\

Let $\mathbf{M}^{\textrm{pre}}_{g}$ (respectively, $\mathbf{M}^{\textrm{pre}}_{g,n}$) be the stack of (resp. $n$-pointed) pre-stable curves of genus $g$, and $\mathcal{C}^{\textrm{pre}}_{g}\longrightarrow \mathbf{M}^{\textrm{pre}}_{g}$ (resp. $\mathcal{C}^{\textrm{pre}}_{g,n}\longrightarrow \mathbf{M}^{\textrm{pre}}_{g,n}$) its universal family (see e.g. \cite{be, op1}).

\begin{df} \begin{itemize} 
\item The \emph{derived stack} $\mathbb{R}\mathbf{M}^{\textrm{pre}}_{g}(X)$ \emph{of prestable maps of genus} $g$ \emph{to} $X$ is defined as $$\mathbb{R}\mathbf{M}^{\textrm{pre}}_{g}(X):= \mathbb{R}\mathrm{HOM}_{\dst /\mathbf{M}^{\textrm{pre}}_{g}} (\mathcal{C}^{\textrm{pre}}_{g}, X\times \mathbf{M}^{\textrm{pre}}_{g}).$$
$\mathbb{R}\mathbf{M}^{\textrm{pre}}_{g}(X)$ is then canonically a derived stack over $\mathbf{M}^{\textrm{pre}}_{g}$, and the corresponding \emph{derived universal family} $\mathbb{R}\mathcal{C}^{\textrm{pre}}_{g;\,X}$ is defined by the following homotopy cartesian square $$\xymatrix{ \mathbb{R}\mathcal{C}^{\textrm{pre}}_{g;\,X}\ar[d] \ar[r] & \mathbb{R}\mathbf{M}^{\textrm{pre}}_{g}(X) \ar[d] \\
\mathcal{C}^{\textrm{pre}}_{g} \ar[r] & \mathbf{M}^{\textrm{pre}}_{g}}$$
\item The \emph{derived stack} $\mathbb{R}\mathbf{M}^{\textrm{pre}}_{g,n}(X)$ \emph{of $n$-pointed prestable maps of genus} $g$ \emph{to} $X$ is defined as $$\mathbb{R}\mathbf{M}^{\textrm{pre}}_{g,n}(X):= \mathbb{R}\mathrm{HOM}_{\dst /\mathbf{M}^{\textrm{pre}}_{g,n}} (\mathcal{C}^{\textrm{pre}}_{g,n}, X\times \mathbf{M}^{\textrm{pre}}_{g,n}).$$
$\mathbb{R}\mathbf{M}^{\textrm{pre}}_{g,n}(X)$ is then canonically a derived stack over $\mathbf{M}^{\textrm{pre}}_{g,n}$, and the corresponding \emph{derived universal family} $\mathbb{R}\mathcal{C}^{\textrm{pre}}_{g,n;\,X}$ is defined by the following homotopy cartesian square $$\xymatrix{ \mathbb{R}\mathcal{C}^{\textrm{pre}}_{g,n;\,X}\ar[d] \ar[r] & \mathbb{R}\mathbf{M}^{\textrm{pre}}_{g,n}(X) \ar[d] \\
\mathcal{C}^{\textrm{pre}}_{g,n} \ar[r] & \mathbf{M}^{\textrm{pre}}_{g,n}}$$
\end{itemize}
\end{df}

\noindent Note that, by definition, $\mathbb{R}\mathcal{C}^{\textrm{pre}}_{g;\,X}$ comes also equipped with a canonical map $$\mathbb{R}\mathcal{C}^{\textrm{pre}}_{g;\,X}\longrightarrow  \mathbb{R}\mathbf{M}^{\textrm{pre}}_{g}(X) \times X .$$ We also have $\tr (\mathbb{R}\mathbf{M}^{\textrm{pre}}_{g}(X)) \simeq \mathbf{M}^{\textrm{pre}}_{g}(X)$ (the stack of prestable maps of genus $g$ to $X$), and  $\tr (\mathbb{R}\mathcal{C}^{\textrm{pre}}_{g;\,X}) \simeq \mathcal{C}^{\textrm{pre}}_{g;\,X}$ (the universal family over the stack of pre-stable maps of genus $g$ to $X$), since the truncation functor $\tr$ commutes with homotopy fibered products. The same is true for the pointed version. \\

We can now use Proposition \ref{opens} to define the derived stable versions. Let $\overline{\mathbf{M}}_{g}(X)$ (respectively, $\overline{\mathbf{M}}_{g,n}(X)$ ) be the open substack of $\mathbf{M}^{\textrm{pre}}_{g}(X)$ (resp. of $\mathbf{M}^{\textrm{pre}}_{g,n}(X)$) consisting of \emph{stable maps of genus $g$ to $X$} (resp. \emph{$n$-pointed stable maps of genus $g$ to $X$}), and $\mathcal{C}_{g;\, X} \longrightarrow \mathbf{M}^{\textrm{pre}}_{g}(X)$ (resp. $\mathcal{C}_{g,n;\, X} \longrightarrow \mathbf{M}^{\textrm{pre}}_{g,n}(X)$) the (induced) universal family (\cite{be, op1}).\\

\begin{df}\label{derstable} 
\begin{itemize}
\item The \emph{derived stack} $\mathbb{R}\overline{\mathbf{M}}_{g}(X)$ \emph{of stable maps of genus} $g$ \emph{to} $X$ is defined as $$\mathbb{R}\overline{\mathbf{M}}_{g}(X):= \phi_{\mathbb{R}\mathbf{M}^{\textrm{pre}}_{g}(X)} (\overline{\mathbf{M}}_{g}(X)).$$
The \emph{derived stable universal family} $$\mathbb{R}\mathcal{C}_{g;\,X}\longrightarrow \mathbb{R}\overline{\mathbf{M}}_{g}(X)$$ is the derived restriction of $\mathbb{R}\mathcal{C}^{\textrm{pre}}_{g;\,X} \rightarrow \mathbb{R}\mathbf{M}^{\textrm{pre}}_{g}(X)$ to  $\mathbb{R}\overline{\mathbf{M}}_{g}(X)$.
\item The \emph{derived stack} $\mathbb{R}\overline{\mathbf{M}}_{g,n}(X)$ \emph{of $n$-pointed stable maps of genus} $g$ \emph{to} $X$ is defined as $$\mathbb{R}\overline{\mathbf{M}}_{g,n}(X):= \phi_{\mathbb{R}\mathbf{M}^{\textrm{pre}}_{g,n}(X)} (\overline{\mathbf{M}}_{g,n}(X)).$$
The \emph{derived stable universal family} $$\mathbb{R}\mathcal{C}_{g,n;\,X}\longrightarrow \mathbb{R}\overline{\mathbf{M}}_{g,n}(X)$$ is the derived restriction of $\mathbb{R}\mathcal{C}^{\textrm{pre}}_{g,n;\,X} \rightarrow \mathbb{R}\mathbf{M}^{\textrm{pre}}_{g,n}(X)$ to  $\mathbb{R}\overline{\mathbf{M}}_{g,n}(X)$.
\end{itemize}
\end{df}

\noindent Recall that 
\begin{itemize}
\item $\tr (\mathbb{R}\overline{\mathbf{M}}_{g}(X)) \simeq \overline{\mathbf{M}}_{g}(X)$;
\item $\tr (\mathbb{R}\mathcal{C}_{g;\,X}) \simeq  \mathcal{C}_{g;\, X}$ ;
\item $\mathbb{R}\mathcal{C}_{g;\,X}$ comes equipped with a canonical map $$\pi: \mathbb{R}\mathcal{C}_{g;\,X}\longrightarrow \mathbb{R}\overline{\mathbf{M}}_{g}(X)\times X ; $$
\item we have a homotopy cartesian diagram in $\dst$ $$\xymatrix{ \overline{\mathbf{M}}_{g}(X) \ar@{^{(}->}[r] \ar[d] & \mathbf{M}^{\textrm{pre}}_{g}(X) \ar[d] \\ \mathbb{R}\overline{\mathbf{M}}_{g}(X) \ar@{^{(}->}[r] & \mathbb{R}\mathbf{M}^{\textrm{pre}}_{g}(X)}$$
\end{itemize}
With the obvious changes, this applies to the pointed version too. \\

Let $g$ a non-negative integer, $\beta \in H_2(X,\mathbb{Z})$, and $\overline{\mathbf{M}}_{g}(X;\beta)$ (respectively, $\overline{\mathbf{M}}_{g,n}(X;\beta)$) be the stack of stable maps (resp. of $n$-pointed stable maps) of type $(g,\beta)$ to $X$ (see e.g. \cite{be} or \cite{op1}); its derived version is given by the following

\begin{df} \begin{itemize}
\item The \emph{derived stack of stable maps of type} $(g,\beta)$ \emph{to} $X$ is defined as the open substack of $\mathbb{R}\overline{\mathbf{M}}_{g}(X)$ $$\mathbb{R}\overline{\mathbf{M}}_{g}(X;\beta):=\phi_{\mathbb{R}\overline{\mathbf{M}}_{g}(X)}(\overline{\mathbf{M}}_{g}(X;\beta)).$$
The \emph{derived stable universal family of type} $(g;\beta)$, $$\mathbb{R}\mathcal{C}_{g, \beta ;\,X}\longrightarrow \mathbb{R}\overline{\mathbf{M}}_{g}(X;\beta) ,$$ is the (derived) restriction of $\mathbb{R}\mathcal{C}_{g;\,X}\longrightarrow \mathbb{R}\overline{\mathbf{M}}_{g}(X)$ to  $\mathbb{R}\overline{\mathbf{M}}_{g}(X;\beta)$.
\item The \emph{derived stack of $n$-pointed stable maps of type} $(g,\beta)$ \emph{to} $X$ is defined as the open substack of $\mathbb{R}\overline{\mathbf{M}}_{g,n}(X)$ $$\mathbb{R}\overline{\mathbf{M}}_{g,n}(X;\beta):=\phi_{\mathbb{R}\overline{\mathbf{M}}_{g,n}(X)}(\overline{\mathbf{M}}_{g,n}(X;\beta)).$$
The \emph{derived stable universal family of type} $(g;\beta)$, $$\mathbb{R}\mathcal{C}_{g,n, \beta ;\,X}\longrightarrow \mathbb{R}\overline{\mathbf{M}}_{g,n}(X;\beta) ,$$ is the (derived) restriction of $\mathbb{R}\mathcal{C}_{g,n;\,X}\longrightarrow \mathbb{R}\overline{\mathbf{M}}_{g,n}(X)$ to  $\mathbb{R}\overline{\mathbf{M}}_{g,n}(X;\beta)$.
\end{itemize}
\end{df}

\noindent Note that, by definition, $\tr (\mathbb{R}\overline{\mathbf{M}}_{g}(X;\beta))\simeq \overline{\mathbf{M}}_{g}(X;\beta),$ therefore $\mathbb{R}\overline{\mathbf{M}}_{g}(X;\beta)$ is a proper derived Deligne-Mumford stack (\cite[2.2.4]{hagII}). Moreover, the derived stable universal family $\mathbb{R}\mathcal{C}_{g, \beta ;\,X}$ comes, by restriction, equipped with a natural map $$\pi: \mathbb{R}\mathcal{C}_{g, \beta ;\,X}\longrightarrow \mathbb{R}\overline{\mathbf{M}}_{g}(X;\beta) \times X .$$ 

We have a homotopy cartesian diagram in $\dst$ $$\xymatrix{ \overline{\mathbf{M}}_{g}(X;\beta) \ar@{^{(}->}[r] \ar[d] & \overline{\mathbf{M}}_{g}(X) \ar[d] \\ \mathbb{R}\overline{\mathbf{M}}_{g}(X;\beta) \ar@{^{(}->}[r] & \mathbb{R}\overline{\mathbf{M}}_{g}(X)}.$$
Analogous remarks are valid in the pointed case.\\

The tangent complex of $\mathbb{R}\overline{\mathbf{M}}_{g}(X;\beta)$ at a stable map $(f:C\rightarrow X)$ of type $(g,\beta)$ (corresponding to a classical point $x_{f}:\mathrm{Spec}(\mathbb{C})\rightarrow \mathbb{R}\overline{\mathbf{M}}_{g}(X;\beta)$) is given by\footnote{The [1] shift in \cite[Thm. 5.4.8]{kcf} is clearly a typo: their proof is correct and yields no shift.} $$\mathbb{T}_{(f:C\rightarrow X)}\simeq \mathbb{R}\Gamma(C, \mathrm{Cone}(\mathbb{T}_{C}\rightarrow f^{*}\mathrm{T}_{X})),$$ where $\mathbb{T}_{C}$ is the tangent complex of $C$ and $T_{X}$ is the tangent sheaf of $X$.\\ The canonical map  $\mathbb{R}\overline{\mathbf{M}}_{g}(X;\beta) \rightarrow \mathbf{M}^{\textrm{pre}}_{g}$ is quasi-smooth. In fact, the fiber at a geometric point, corresponding to prestable curve $C$, is the derived stack $\mathbb{R}\mathrm{HOM}_{\beta}(C,X)$ whose tangent complex at a point $f:C\rightarrow X$ is $\mathbb{R}\Gamma(C,f^*T_S)$ which, obviously, has cohomology only in degrees $[0,1]$. But $\mathbf{M}^{\textrm{pre}}_{g}$ is smooth, and any  derived stack quasi-smooth over a smooth base is quasi-smooth (by the corresponding exact triangle of tangent complexes). Therefore the derived stack $\mathbb{R}\overline{\mathbf{M}}_{g}(X;\beta)$ is \emph{quasi-smooth}. \\
Proposition \ref{inducedobstruction} then recovers the \emph{standard} (absolute) perfect obstruction theory on $\overline{\mathbf{M}}_{g}(X;\beta)$ via the canonical map $$j^{*} (\mathbb{L}_{\mathbb{R}\overline{\mathbf{M}}_{g}(X;\beta)}) \longrightarrow \mathbb{L}_{\overline{\mathbf{M}}_{g}(X;\beta)}$$ induced by the closed immersion $j: \overline{\mathbf{M}}_{g}(X;\beta)\hookrightarrow \mathbb{R}\overline{\mathbf{M}}_{g}(X;\beta)$. \\ 

In the pointed case, the tangent complex of $\mathbb{R}\overline{\mathbf{M}}_{g,n}(X;\beta)$ at a pointed stable map $(f:(C ; x_1,\ldots, x_n)\rightarrow X)$ of type $(g,\beta)$ (corresponding to a classical point $x_{f}:\mathrm{Spec}(\mathbb{C})\rightarrow \mathbb{R}\overline{\mathbf{M}}_{g,n}(X;\beta)$) is likewise given by $$\mathbb{T}_{(f:(C ; x_1,\ldots, x_n)\rightarrow X)}\simeq \mathbb{R}\Gamma(C, \mathrm{Cone}(\mathbb{T}_{C}(-\sum_{i} x_{i})\rightarrow f^{*}\mathrm{T}_{X})).$$ The same argument as above proves that also the canonical map  $\mathbb{R}\overline{\mathbf{M}}_{g,n}(X;\beta) \rightarrow \mathbf{M}^{\textrm{pre}}_{g,n}$ is quasi-smooth, and Proposition \ref{inducedobstruction} then recovers the \emph{standard} absolute perfect obstruction theory on $\overline{\mathbf{M}}_{g,n}(X;\beta)$ via the canonical map $$j^{*} (\mathbb{L}_{\mathbb{R}\overline{\mathbf{M}}_{g,n}(X;\beta)}) \longrightarrow \mathbb{L}_{\overline{\mathbf{M}}_{g,n}(X;\beta)}$$ induced by the closed immersion $j: \overline{\mathbf{M}}_{g,n}(X;\beta)\hookrightarrow \mathbb{R}\overline{\mathbf{M}}_{g,n}(X;\beta)$. Note that, as observed in \cite[5.3.5]{op1}, this obstruction theory yields trivial Gromov-Witten invariants on $\overline{\mathbf{M}}_{g,n}(X;\beta)$ for $X=S$ a $K3$ surface. Hence the need for another obstruction theory carrying more interesting curves counting invariants on a $K3$-surface: this will be the so-called \emph{reduced} obstruction theory (see \S \ref{review}, \S \ref{our}, and Theorem \ref{quasismooth}). \\

\section{The derived determinant morphism}
In this section we start by defining a quite general \emph{perfect determinant map} of derived stacks $$\textrm{det}_{\textrm{Perf}}: \mathbf{Perf}\longrightarrow \mathbf{Pic} = B\mathbb{G}_{m}$$ whose construction requires a small detour into Waldhausen $K$-theory. We think this perfect determinant might play an important role in other contexts as well, e.g. in a general GW/DT correspondence.

Using the perfect determinant together with a natural perfect complex on $\mathbb{R}\overline{\mathbf{M}}_{g}(X;\beta)$, we will be able to define a map  $$\delta_{1}(X): \mathbb{R}\overline{\mathbf{M}}_{g}(X)\longrightarrow \mathbb{R}\mathbf{Pic}(X)$$ which will be one of the main ingredients in the construction of the \emph{reduced} derived stack of stable maps $\mathbb{R}\overline{\mathbf{M}}^{\textrm{red}}_{g}(S;\beta)$, for a $K3$-surface $S$, given in the next section. \\

\subsection{The perfect determinant map}
The aim of this subsection is to produce a determinant morphism $\det_{\mathrm{Perf}} : \mathbf{Perf}\longrightarrow \mathbf{Pic}$ in $\mathrm{Ho}(\dst)$ extending the natural determinant morphism $\mathbf{Vect}\longrightarrow \mathbf{Pic}.$ To do this, we will have to pass through Waldhausen $K$-theory.\\

We start with the classical determinant map in $\mathrm{Ho}(\st)$, $\det : \mathbf{Vect}\longrightarrow \mathbf{Pic}$, induced by the map sending a vector bundle to its top exterior power. Consider the following simplicial stacks $$B_{\bullet}\mathbf{Pic}: \Delta^{\textrm{op}} \ni [n] \longmapsto (\mathbf{Pic})^{n}$$ (with the simplicial structure maps given by tensor products of line bundles, or equivalently, induced by the product in the group structure of $B\mathbb{G}_{m}\simeq \mathbf{Pic}$), and 
$$B_{\bullet}\mathbf{Vect}: \Delta^{\textrm{op}} \ni [n] \longmapsto wS_{n}\mathbf{Vect},$$ where, for any commutative $\mathbb{C}$-algebra $R$, $wS_{n}\mathbf{Vect} (R)$ is the nerve of the category of sequences of split monomorphisms $$0\rightarrow M_{1} \rightarrow M_{2} \rightarrow \ldots \rightarrow M_{n}\rightarrow 0$$ with morphisms the obvious equivalences, and the simplicial structure maps are the natural ones described in \cite[1.3]{wal}. Similarly, we define the simplicial object in stacks $$B_{\bullet}\mathbf{Perf}: \Delta^{\textrm{op}} \ni [n] \longmapsto wS_{n}\mathbf{Perf}$$ (see \cite[1.3]{wal} for the definition of $wS_{n}$ in this case). Now, $B_{\bullet}\mathbf{Pic}$ and  $B_{\bullet}\mathbf{Vect}$, and $B_{\bullet}\mathbf{Perf}$ are pre-$\Delta^{\textrm{op}}$-stacks according to Def. 1.4.1 of \cite{chaff}, and  the map $\det$ extends to a morphism $$\textrm{det}_{\bullet}: B_{\bullet}\mathbf{Vect}\longrightarrow B_{\bullet}\mathbf{Pic}$$ in the homotopy category of pre-$\Delta^{\textrm{op}}$-stacks. By applying the functor $i:\mathrm{Ho}(\st)\rightarrow \mathrm{Ho}(\dst)$ (that will be, according to our conventions, omitted from notations), we get a determinant morphism (denoted in the same way)  $$\textrm{det}_{\bullet} : B_{\bullet}\mathbf{Vect}\longrightarrow B_{\bullet}\mathbf{Pic}$$ in the homotopy category of pre-$\Delta^{\textrm{op}}$-derived stacks. We now pass to Waldhausen $K-$theory, i.e. apply $\textrm{K}:=\Omega \circ | - |$ (see \cite[Thm 1.4.3]{chaff}, where the loop functor $\Omega$ is denoted by $\mathbb{R}\Omega_{*}$, and the realization functor $| - |$ by $B$), and observe that, by \cite[Thm 1.4.3 (2)]{chaff}, there is a canonical isomorphism in $\mathrm{Ho}(\dst)$ $$\textrm{K}(B_{\bullet}\mathbf{Pic}) \simeq \mathbf{Pic}$$ since $\mathbf{Pic}$ is group-like (i.e. an $H_{\infty}$-stack in the parlance of \cite[Thm 1.4.3]{chaff}). This gives us a map in $\mathrm{Ho}(\dst)$ 
$$\mathrm{K}(\textrm{det}_{\bullet}):\textrm{K}(B_{\bullet}\mathbf{Vect}) \longrightarrow \mathbf{Pic}. $$ Now, consider the map $u: \mathbf{K}^{\textrm{Vect}}:=\textrm{K}(B_{\bullet}\mathbf{Vect}) \longrightarrow \textrm{K}(B_{\bullet}\mathbf{Perf}):=\mathbf{K}^{\textrm{Perf}}$  in $\mathrm{Ho}(\dst)$, induced by the inclusion $\mathbf{Vect}\hookrightarrow \mathbf{Perf}$. By \cite[Thm. 1.7.1]{wal}, $u$ is an isomorphism in $\mathrm{Ho}(\dst)$. Therefore, we get a diagram in $\mathrm{Ho}(\dst)$ \\

$$\xymatrix{ & \mathbf{K}^{\textrm{Vect}} \ar[r]^-{\mathrm{K}(\textrm{det}_{\bullet})} \ar[d]_-{u} & \mathbf{Pic} \\
\mathbf{Perf} \ar[r]_-{\textrm{1st-level}} & \mathbf{K}^{\textrm{Perf}} &}
$$ \\ 

\noindent where $u$ is an isomorphism. This allows us to give the following

\begin{df}\label{det} The induced map in $\mathrm{Ho}(\dst)$ $$\mathrm{det}_{\mathrm{Perf}} : \mathbf{Perf}\longrightarrow \mathbf{Pic}$$ is called the \emph{perfect determinant} morphism.
\end{df}

For any complex scheme (or derived stack) $X$, the perfect determinant morphism $\det_{\mathrm{Perf}} : \mathbf{Perf}\longrightarrow \mathbf{Pic}$ induces a map in  $\mathrm{Ho}(\dst)$ $$\mathrm{det}_{\mathrm{Perf}}(X): \mathbb{R}\mathbf{Perf}(X):=\mathbb{R}\mathrm{HOM}_{\dst}(X, \mathbf{Perf})\longrightarrow \mathbb{R}\mathrm{HOM}_{\dst}(X, \mathbf{Pic})=:\mathbb{R}\mathbf{Pic}(X).$$
As perhaps not totally unexpected (e.g. \cite[Rem. 5.3.3]{ill}), the tangent morphism to the perfect determinant map is given by the \emph{trace for perfect complexes}

\begin{prop} \label{trace} Let $X$ be a complex quasi-projective scheme, and $\mathrm{det}_{\mathrm{Perf}}(X) : \mathbb{R}\mathbf{Perf}(X) \rightarrow \mathbb{R}\mathbf{Pic}(X)$ the induced perfect determinant map. For any complex point $x_{E}: \mathrm{Spec}\, \mathbb{C} \rightarrow \mathbb{R}\mathbf{Perf}(X)$, corresponding to a perfect complex $E$ over $X$, the tangent map $$\mathbb{T}_{x_{E}}\mathrm{det}_{\mathrm{Perf}}(X) : \mathbb{T}_{x_{E}}\mathbb{R}\mathbf{Perf}(X) \simeq \mathbb{R}\mathrm{Hom}(E,E)[1] \longrightarrow \mathbb{R}\mathrm{Hom}(\mathcal{O}_{S},\mathcal{O}_{S})[1]\simeq \mathbb{T}_{x_{E}}\mathbb{R}\mathbf{Pic}(X)$$ is given by \emph{$\mathrm{tr}_{E}[1]$}, where \emph{$\mathrm{tr}_{E}$} is the trace map for the perfect complex $E$ of \emph{\cite[Ch. 5, 3.7.3]{ill}.}
\end{prop}
\noindent\textbf{Proof.} Let $\mathbb{R}\mathbf{Perf}^{\textrm{strict}}(X):=\mathbb{R}\mathrm{HOM}_{\dst}(X,\mathbf{Perf}^{\textrm{strict}})$ be the derived stack of \emph{strict} perfect complexes on $X$ (\cite[Exp. I, 2.1]{sga6}). Since $X$ is quasi-projective, the canonical map $\mathbb{R}\mathbf{Perf}^{\textrm{strict}}(X) \rightarrow \mathbb{R}\mathbf{Perf}(X)$ is an isomorphism in $\mathrm{Ho}(\dst)$. Therefore (e.g. \cite[Exp. I, 8.1.2]{sga6}), the comparison statement is reduced to the case where $E$ is a vector bundle on $X$, which is a direct computation and is left to the reader.
\hfill $\Box$ \\

\subsection{The map $\mathbb{R}\overline{\mathbf{M}}_{g}(X)\longrightarrow \mathbb{R}\mathbf{Perf}(X)$}\label{AX}

A map $$\mathbb{R}\overline{\mathbf{M}}_{g}(X)\longrightarrow \mathbb{R}\mathbf{Perf}(X) = \mathbb{R}\mathrm{HOM}_{\dst}(X, \mathbf{Perf})$$ in  $\mathrm{Ho}(\dst)$ is, by adjunction, the same thing as a map $$\mathbb{R}\overline{\mathbf{M}}_{g}(X)\times X\longrightarrow \mathbf{Perf}$$ i.e. a perfect complex on $\mathbb{R}\overline{\mathbf{M}}_{g}(X)\times X$; so, it is enough to find such an appropriate perfect complex.\\

Let $$\pi:\mathbb{R}\mathcal{C}_{g ;\,X}\longrightarrow \mathbb{R}\overline{\mathbf{M}}_{g}(X) \times X$$ be the derived stable universal family (\S \ref{derivedstable}). 

\begin{prop} \label{perf} $\mathbb{R}\pi_{*}(\mathcal{O}_{\mathbb{R}\mathcal{C}_{g;\,X}})$ is a perfect complex on $\mathbb{R}\overline{\mathbf{M}}_{g}(X)\times X$. 
\end{prop}

\noindent \textbf{Proof.} The truncation of $\pi$ is proper, hence it is enough to prove that $\pi$ is also quasi-smooth. To see this, observe that both $\mathbb{R}\mathcal{C}_{g;\,X}$ and $\mathbb{R}\overline{\mathbf{M}}_{g}(X) \times X$ are smooth over $\mathbb{R}\overline{\mathbf{M}}_{g}(X)$. Then we conclude, since any map between derived stacks smooth over a base is quasi-smooth.

\hfill $\Box$ \\

\begin{rmk} \emph{If we fix a class $\beta \in H_{2}(X,\mathbb{Z})$, the corresponding $\beta$-decorated version of Proposition \ref{perf} obviously holds.}
\end{rmk}

We may therefore give the following

\begin{df} We will denote by $$ \mathrm{A}_{X}:  \mathbb{R}\overline{\mathbf{M}}_{g}(X)\longrightarrow \mathbb{R}\mathbf{Perf}(X)$$ the map induced by the perfect complex $\mathbb{R}\pi_{*}(\mathcal{O}_{\mathbb{R}\mathcal{C}_{g;\,X}})$.
\end{df}

Note that, by definition, $\mathrm{A}_{X}$ sends a complex point of  $\mathbb{R}\overline{\mathbf{M}}_{g}(X)$, corresponding to a stable map $f:C\rightarrow X$ to the perfect complex $\mathbb{R}f_{*}\mathcal{O}_{C}$ on $X$.\\

\noindent \textbf{The tangent morphism of  $\mathrm{A}_{X}$.} The tangent morphism of $\mathrm{A}_{X}$ is related to the Atiyah class of $\mathbb{R}\pi_{*}(\mathcal{O}_{\mathbb{R}\mathcal{C}_{g;\,X}})$, and pointwise on $\mathbb{R}\overline{\mathbf{M}}_{g}(X)$ to the Atiyah class map of the perfect complex $\mathbb{R}f_{*}\mathcal{O}_{C}$: this is explained in detail in Appendix \ref{AXexplained}, so we will just recall here the results and the notations we will need in the rest of the main text. \\
Let us write $\mathcal{E}:=\mathbb{R}\pi_{*}(\mathcal{O}_{\mathbb{R}\mathcal{C}_{g;\,X}})$; since this is a perfect complex on $\mathbb{R}\overline{\mathbf{M}}_{g}(X)\times X$, its Atiyah class map (see Appendix \ref{AXexplained}) $$\mathrm{at}_{\mathcal{E}}: \mathcal{E}\longrightarrow \mathbb{L}_{\mathbb{R}\overline{\mathbf{M}}_{g}(X)\times X}\otimes \mathcal{E}[1]$$ corresponds uniquely, by adjunction, to a map, denoted in the same way, $$\mathrm{at}_{\mathcal{E}}: \mathbb{T}_{\mathbb{R}\overline{\mathbf{M}}_{g}(X)\times X}\longrightarrow \mathcal{E}^{\vee}\otimes \mathcal{E}[1].$$ Let $x$ be a complex point $x$ of $\mathbb{R}\overline{\mathbf{M}}_{g}(X)$ corresponding to a stable map $f:C\rightarrow X$, and let $p: C\rightarrow \mathrm{Spec}\, \mathbb{C}$ and $q:X \rightarrow \mathrm{Spec}\, \mathbb{C}$ denote the structural morphisms, so that $p=q\circ f$. Correspondingly, we have a ladder of homotopy cartesian diagrams 
$$\xymatrix{C \ar[r]^-{\iota_{f}} \ar[d]_{f} & \mathbb{R}\mathcal{C}_{g;\,X} \ar[d]^-{\pi}  & \\ X \ar[r]^-{\underline{x}} \ar[d]_-{q} & \mathbb{R}\overline{\mathbf{M}}_{g}(X)\times X \ar[d]^-{\mathrm{pr}} \ar[r]^-{\mathrm{pr}_{X}} & X \ar[d]^-{q} \\ \mathrm{Spec}\, \mathbb{C} \ar[r]_-{x} & \mathbb{R}\overline{\mathbf{M}}_{g}(X) \ar[r] & \mathrm{Spec}\, \mathbb{C}}$$ Let us consider the perfect complex $E:=\mathbb{R}f_{*}\mathcal{O}_{C}$ on $X$. By \cite[Ch. 4, 2.3.7]{ill}, the complex $E$ has an Atiyah class map 
$$\mathrm{at}_{E}: E \longrightarrow E\otimes \Omega_{X}^{1}[1]$$ which corresponds uniquely ($E$ being perfect) by adjunction to a map (denoted in the same way) $$\mathrm{at}_{E}: T_{X}\longrightarrow \mathbb{R}\underline{\mathrm{End}}_{X}(\mathbb{R}f_{*}\mathcal{O}_{C})[1].$$

\begin{prop}\label{finallyatiyah} In the situation and notations above, we have that 
\begin{itemize}
\item the tangent map of \emph{$\mathrm{A}_{X}$} fits into the following commutative diagram $$\xymatrix{\mathbb{T}_{\mathbb{R}\overline{\mathbf{M}}_{g}(X)} \ar[r]^-{\mathbb{T}\mathrm{A}_{X}} \ar[d]_-{\textrm{can}} & \mathrm{A}_{X}^{*}\mathbb{T}_{\mathbb{R}\mathbf{Perf}(X)} \ar[r]^-{\sim} & \mathbb{R}pr_{*}(\mathcal{E}^{\vee}\otimes \mathcal{E})[1] \\ \mathbb{R}pr_{*}pr^{*}\mathbb{T}_{\mathbb{R}\overline{\mathbf{M}}_{g}(X)} \ar[r]_-{\textrm{can}} & \mathbb{R}pr_{*}(pr^{*}\mathbb{T}_{\mathbb{R}\overline{\mathbf{M}}_{g}(X)}\oplus pr_{X}^{*}T_{X})  \ar[r]_-{\sim} & \mathbb{R}pr_{*}\mathbb{T}_{\mathbb{R}\overline{\mathbf{M}}_{g}(X) \times X} \ar[u]_-{\mathbb{R}pr_{*}(\mathrm{at}_{\mathcal{E})} }}
$$ where \emph{can} denote obvious canonical maps, and $\mathcal{E}:=\mathbb{R}\pi_{*}(\mathcal{O}_{\mathbb{R}\mathcal{C}_{g;\,X}})$.
\item The tangent map to $\mathrm{A}_{X}$ at $x=(f:C\rightarrow X)$, is the composition $$\xymatrix{\mathbb{T}_{x}\mathrm{A}_{X}: \mathbb{T}_{x} \mathbb{R}\overline{\mathbf{M}}_{g}(X) \simeq \mathbb{R}\Gamma(C, \mathrm{Cone}(\mathbb{T}_{C}\rightarrow f^{*}T_{X})) \ar[r] & \mathbb{R}\Gamma(X, \underline{x}^{*}\mathbb{T}_{\mathbb{R}\overline{\mathbf{M}}_{g}(X) \times X}) \ar[r]  & }$$ $$\xymatrix{ \ar[rrr]^-{\mathbb{R}\Gamma(X,\underline{x}^{*}\textrm{at}_{E})} & & & \mathbb{R}\mathrm{End}_{X}(\mathbb{R}f_{*}\mathcal{O}_{C})[1] \simeq \mathbb{T}_{\mathbb{R}f_*\mathcal{O}_{C}}\mathbb{R}\mathbf{Perf}(X)}$$ where $E:=\mathbb{R}f_{*}\mathcal{O}_{C}$  
\item The composition $$\xymatrix{\mathbb{R}\Gamma(X,T_X) \ar[r]^-{\textrm{can}} & \mathbb{R}\Gamma(X,\mathbb{R}f_{*}f^{*}T_X) \ar[r]^-{\textrm{can}} &   \mathbb{R}\Gamma(X,\textrm{Cone}(\mathbb{R}f_{*}\mathbb{T}_{C}\rightarrow \mathbb{R}f_{*}f^{*}T_X))\simeq \mathbb{T}_{x} \mathbb{R}\overline{\mathbf{M}}_{g}(X) \ar[r] & }$$ $$\xymatrix{ \ar[r]^-{\mathbb{T}_{x}\mathrm{A}_{X}} & x^{*}\mathrm{A}_{X}^{*}\mathbb{T}\mathbb{R}\mathbf{Perf}(X)\simeq \mathbb{T}_{\mathbb{R}f_{*}\mathcal{O}_{C}}\mathbb{R}\mathbf{Perf}(X) \simeq \mathbb{R}\mathrm{End}_{X}(\mathbb{R}f_{*}\mathcal{O}_{C})[1]  }$$ coincides with $\mathbb{R}\Gamma(X, \textrm{at}_{E})$, where $E:=\mathbb{R}f_{*}\mathcal{O}_{C}$.
\end{itemize}
\end{prop}
\noindent \textbf{Proof.} See Appendix \ref{AXexplained}.
\hfill $\Box$ \\ 

\begin{df} \label{delta1} We denote by $\delta_{1}(X)$ the composition
 $$\xymatrix{\mathbb{R}\overline{\mathbf{M}}_{g}(X)\ar[r]^-{\mathrm{A}_{X}} &\mathbb{R}\mathbf{Perf}(X) \ar[rr]^-{\mathrm{det}_{\mathrm{Perf}}(X)} & & \mathbb{R}\mathbf{Pic}(X),}$$ and, for a complex point $x$ of $\mathbb{R}\overline{\mathbf{M}}_{g}(X)$ corresponding to a stable map $f:C\rightarrow X$, by $$\xymatrix{\Theta_{f}:= \mathbb{T}_{f}\delta_{1}(X): \mathbb{T}_{(f:C \rightarrow X)} \mathbb{R}\overline{\mathbf{M}}_{g}(X) \ar[r]^-{\mathbb{T}_{x}\mathrm{A}_{X}} & \mathbb{T}_{\mathbb{R}f_{*}\mathcal{O}_{C}}\mathbb{R}\mathbf{Perf}(X) \ar[r]^-{\mathrm{tr}_{X}} & \mathbb{T}_{\det (\mathbb{R}f_{*}\mathcal{O}_{C})} \mathbb{R}\mathbf{Pic}(X).}$$ 
\end{df}

\bigskip

Note that, as a map of explicit complexes, we have $$\xymatrix{\Theta_{f} : \mathbb{R}\Gamma(C,\mathrm{Cone}(\mathbb{T}_{C}\rightarrow f^{*}T_{X})) \ar[r]^-{\mathbb{T}_{x}\mathrm{A}_{X}} & \mathbb{R}\mathrm{Hom}_{X}(\mathbb{R}f_{*}\mathcal{O}_{C}, \mathbb{R}f_{*}\mathcal{O}_{C})[1] \ar[r]^-{\mathrm{tr}_{X}} & \mathbb{R}\Gamma (X, \mathcal{O}_{X})[1]}$$ \\

\begin{rmk} \label{atiyah2}\emph{\textbf{- First Chern class of $\mathbb{R}f_{*}\mathcal{O}_{C}$ and the map $\Theta_{f}$.} Using Proposition \ref{finallyatiyah}, we can relate the map $\Theta_{f}$ above to the \emph{first Chern class} of the perfect complex $\mathbb{R}f_{*}\mathcal{O}_{C}$ (\cite[Ch. V]{ill}). With the same notations as in Prop. \ref{finallyatiyah}, the following diagram is commutative $$\xymatrix{\mathbb{R}q_{*}T_{X} \ar[d] \ar[r]^-{\mathbb{R}q_{*}(\mathrm{at}_{\mathbb{R}f_{*}\mathcal{O}_{C}})} & \mathbb{R}q_{*}\mathbb{R}\underline{\mathrm{End}}_{X}(\mathbb{R}f_{*}\mathcal{O}_{C})[1]  \ar[r]^-{\mathrm{tr}} & \mathbb{R}q_{*}\mathcal{O}_{X}[1] \\   \mathbb{R}q_{*} \mathbb{R}f_{*}f^{*}T_{X} \simeq  \mathbb{R}p_{*}f^{*}T_{X} \ar[r] & \mathbb{R}p_{*}\mathrm{Cone}(\mathbb{T}_{C}\rightarrow f^{*}T_{X}) \ar[u]_-{\mathbb{T}_{x}\mathrm{A}_{X}} \ar[r]^-{\Theta_{f}} & \mathbb{R}q_{*}\mathcal{O}_{X}[1] \ar[u]_{\mathrm{id}}}.$$ In this diagram, the composite upper row is the image under $\mathbb{R}q_{*}$ of the first Chern class $c_{1}(\mathbb{R}f_{*}\mathcal{O}_{C}) \in \mathrm{Ext}_{X}^{1}(T_{X}, \mathcal{O}_{X}) \simeq H^{1}(X, \Omega_{X}^{1})$.} 
\end{rmk}

\noindent \textbf{Pointed case -} In the pointed case, if  $$\pi:\mathbb{R}\mathcal{C}_{g,n ;\,X}\longrightarrow \mathbb{R}\overline{\mathbf{M}}_{g,n}(X) \times X$$ is the derived stable universal family (\S \ref{derivedstable}), the same argument as in Proposition \ref{perf} shows that $\mathbb{R}\pi_{*}(\mathcal{O}_{\mathbb{R}\mathcal{C}_{g;\,X}})$ is a perfect complex on $\mathbb{R}\overline{\mathbf{M}}_{g,n}(X)\times X$. And we give the analogous

\begin{df} \label{delta1pointed} 
\begin{itemize}
\item We denote by $$\mathrm{A}^{(n)}_{X}: \mathbb{R}\overline{\mathbf{M}}_{g,n}(X)\longrightarrow \mathbb{R}\mathbf{Perf}(X)$$ the map induced by the perfect complex $\mathbb{R}\pi_{*}(\mathcal{O}_{\mathbb{R}\mathcal{C}_{g,n;\,X}})$.
\item We denote by $\delta^{(n)}_{1}(X)$ the composition
 $$\xymatrix{\mathbb{R}\overline{\mathbf{M}}_{g,n}(X)\ar[r]^-{\mathrm{A}^{(n)}_{X}} &\mathbb{R}\mathbf{Perf}(X) \ar[rr]^-{\mathrm{det}_{\mathrm{Perf}}(X)} & & \mathbb{R}\mathbf{Pic}(X),}$$ and, for a complex point $x$ of $\mathbb{R}\overline{\mathbf{M}}_{g}(X)$ corresponding to a pointed stable map $f:(C; x_1, \ldots, x_n)\rightarrow X$, by $$\xymatrix{\Theta^{(n)}_{f}:= \mathbb{T}_{f}\delta^{(n)}_{1}(X): \mathbb{T}_{f} \mathbb{R}\overline{\mathbf{M}}_{g,n}(X) \ar[r]^-{\mathbb{T}_{x}\mathrm{A}^{(n)}_{X}} & \mathbb{T}_{\mathbb{R}f_{*}\mathcal{O}_{C}}\mathbb{R}\mathbf{Perf}(X) \ar[r]^-{\mathrm{tr}_{X}} & \mathbb{T}_{\det (\mathbb{R}f_{*}\mathcal{O}_{C})} \mathbb{R}\mathbf{Pic}(X).}$$
\end{itemize}
\end{df}

And again, if we fix a class $\beta \in H_{2}(X,\mathbb{Z})$, we have the corresponding $\beta$-decorated version of Definition \ref{delta1pointed}. 

\section{The reduced derived stack of stable maps to a $K3$-surface}

In this section we specialize to the case of a $K3$-surface $S$, with a fixed nonzero curve class $\beta \in H_{2}(S;\mathbb{Z})\simeq H^{2}(S;\mathbb{Z})$. After recalling in some detail the reduced obstruction theory of O-M-P-T, we first identify canonically the derived Picard stack $\mathbb{R}\mathbf{Pic}(S)$  with  $\mathbf{Pic}(S) \times \mathbb{R}\mathrm{Spec}(\mathrm{Sym} (H^{0}(S, K_{S})[1]))$ where $K_{S}$ is the canonical sheaf of $S$. This result is then used to define the \emph{reduced} version $\mathbb{R}\overline{\mathbf{M}}^{\textrm{red}}_{g}(S;\beta)$ of the derived stack of stable maps of type $(g, \beta)$ to $S$ (and its $n$-pointed variant $\mathbb{R}\overline{\mathbf{M}}^{\textrm{red}}_{g,n}(S;\beta)$), and to show that this induces, via the canonical procedure available for any algebraic derived stack, a modified  obstruction theory on its truncation $\overline{\mathbf{M}}_{g}(S;\beta)$ whose deformation and obstruction spaces are then compared with those of the  reduced theory of O-M-P-T. As a terminological remark, given an obstruction theory, we will call \emph{deformation space} what is usually called its \emph{tangent space} (while we keep the terminology \emph{obstruction space}). We do this to avoid confusion with tangent spaces, tangent complexes or tangent cohomologies of possibly related (derived) stacks.\\

\subsection{Review of reduced obstruction theory}\label{review}
For a $K3$-surface $S$, the moduli of stable maps of genus $g$ curves to $S$ with non-zero effective class $\beta \in H^{1,1}(S,\mathbb{C}) \cap H^2(S,\mathbb{Z})$ (note that Poincar\'e duality yields a canonical isomorphism $H_{2}(S;\mathbb{Z})\simeq H^{2}(S;\mathbb{Z})$ between singular (co)homologies) carries a relative perfect obstruction theory. This obstruction theory is given by 
\[
 (R \pi _* F^* T_S)^\vee \to \mathbb{L} _{\overline{\mathbf{M}}_{g}(S;\beta) / \mathbf{M}^{\textrm{pre}}_g}. \ 
\]
Here $\pi \colon \mathcal{C}_{g, \beta ;\,S} \to \overline{\mathbf{M}}_{g}(S;\beta)$ is the universal curve, $F \colon \mathcal{C}_{g, \beta ;\,S} \to S$ is the universal morphism from the universal curve to $S$, and $\mathbf{M}^{\textrm{pre}}_g$ denotes the Artin stack of prestable curves. A Riemann-Roch argument along with the fact that a $K3$-surface has trivial canonical bundle yields the expected dimension of $\overline{\mathbf{M}}_{g}(S;\beta)$:
\[
 \mathrm{exp \, dim} \; \overline{\mathbf{M}}_{g}(S;\beta) = g-1.
\]
We thus expect no rational curves on a $K3$-surface. This result stems from the deformation invariance of Gromov-Witten invariants. A $K3$-surface admits deformations such that the homology class $\beta$ is no longer of type $(1,1)$, and thus can not be the class of a curve. 

This is unfortunate, given the rich literature on enumerative geometry of $K3$-surfaces, and is in stark contrast to the well-known conjecture that a projective $K3$-surface over an algebraically closed field contains infinitely many rational curves. Further evidence that there should be an interesting Gromov-Witten theory of $K3$-surfaces are the results of Bloch, Ran and Voisin that rational curves deform in a family of $K3$-surfaces provided their homology classes remain of type $(1,1)$. The key ingredients in the proof is the semi-regularity map. We thus seek a new kind of obstruction theory for $\overline{\mathbf{M}}_{g}(S;\beta)$ which is deformation invariant only for such deformations of $S$ which keep $\beta$ of type $(1,1)$. 

Such a new obstruction theory, called the \emph{reduced obstruction theory}, was introduced in \cite{op, mp, mpt}. Sticking to the case of moduli of morphisms from a \emph{fixed} curve $C$ to $S$, the obstruction space at a fixed morphism $f$ is $H^1 (C, f^* T_S)$. \\
This obstruction space admits a map

$$\xymatrix{H^1 (C, f^* T_S) \ar[r]^-{\sim} & H^1 (C, f^* \Omega_S) \ar[r]^-{H^{1}(df)} & H^1 (C, \Omega^{1} _{C}) \ar[r] & H^1 (C, \omega _C) \simeq \mathbb{C}},$$ 

\noindent where the first isomorphism is induced by the choice of a holomorphic symplectic form on $S$. The difficult part is to prove that all obstructions for all types of deformations of $f$ (and not only curvilinear ones) lie in the kernel of this map. Once this is proven, $\overline{\mathbf{M}}_{g}(S;\beta)$ carries a reduced obstruction theory which yields a virtual class, called the \emph{reduced class}. This reduced class is one dimension larger that the one obtained from the standard perfect obstruction theory and leads to many interesting enumerative results (see \cite{pak3, mp, mpt}).  \\

We will review below the construction of the reduced deformation and obstruction spaces giving all the details will be needed in our comparison result (Thm. \ref{quasismooth}).

\subsubsection{Deformation and obstruction spaces of the reduced theory according to O-M-P-T} \label{omptspaces}
For further reference, we give here a self-contained treatment of the reduced \emph{deformation} and reduced \emph{obstruction} spaces on $\overline{\mathbf{M}}_{g}(S;\beta)$ according to Okounkov-Maulik-Pandharipande-Thomas.\\

Let us fix a stable map $f:C\rightarrow S$ of class $\beta \neq 0$ and genus $g$; $p: C\rightarrow \mathrm{Spec}\, \mathbb{C}$ and $q:S \rightarrow \mathrm{Spec}\, \mathbb{C}$ will denote the structural morphisms. Let $\underline{\omega}_{C} \simeq p^{!}\mathcal{O}_{\mathrm{Spec}\, \mathbb{C}}$ be the dualizing complex of $C$, and $\omega_{C}=\underline{\omega}_{C}[-1]$ the corresponding dualizing sheaf.\\

First of all, the \emph{deformation spaces} of the standard (i.e. unreduced) and reduced theory, at the stable map $f$, coincide with $$H^{0}(C,\mathrm{Cone}(\mathbb{T}_{C} \rightarrow f^{*}T_{S}))$$ where $\mathbb{T}_{C}$ is the cotangent complex of the curve $C$.\\

Let's recall now (\cite[\S 3.1]{pak3}) the construction of the reduced \emph{obstruction space}. We give here a canonical version, independent of the choice of a holomorphic symplectic form on $S$. \\
Consider the canonical isomorphism\footnote{We use throughout the standard abuse of writing $\mathcal{F}\otimes V$ for $\mathcal{F}\otimes_{\mathcal{O}_{X}} p^{*}V$, for any scheme $p:X\rightarrow \mathrm{Spec}\, \mathbb{C}$, any $\mathcal{O}_{X}$-Module $\mathcal{F}$, and any $\mathbb{C}$-vector space $V$.} $$\xymatrix{\varphi: T_{S} \otimes H^{0}(S,K_{S}) \ar[r]^-{\sim} &\Omega_{S}^{1}.}$$ By tensoring this by $H^{0}(S,K_{S})^{\vee} \simeq H^{2}(S,\mathcal{O}_{S})$ (this isomorphism is canonical by Serre duality) which is of dimension $1$ over $\mathbb{C}$, we get a canonical sequence of isomorphisms of $\mathcal{O}_{S}$-Modules $$\xymatrix{T_{S} & T_{S} \otimes H^{0}(S,K_{S}) \otimes H^{2}(S,\mathcal{O}_{S}) \ar[l]_-{\sim} \ar[r]^-{\sim} & \Omega_{S}^{1} \otimes H^{2}(S,\mathcal{O}_{S}).}$$ We denote by $\psi: T_{S} \rightarrow \Omega_{S}^{1} \otimes H^{2}(S,\mathcal{O}_{S})$ the induced, canonical isomorphism. Form this, we get a canonical  isomorphism of $\mathcal{O}_{C}$-Modules $$\xymatrix{f^{*}\psi : f^{*}T_{S} \ar[r]^-{\sim} & f^{*}(\Omega_{S}^{1}) \otimes H^{2}(S,\mathcal{O}_{S}).}$$ 

Now consider the canonical maps $$\xymatrix{f^{*}\Omega^{1}_{S} \ar[r]^-{s} & \Omega^{1}_{C} \ar[r]^-{t} & \omega_{C}\simeq p^{!}\mathcal{O}_{\mathrm{Spec}\, \mathbb{C}}[-1]}$$ where $\omega_{C}\simeq \underline{\omega}_{C}[-1]$ is the dualizing sheaf of $C$ and $\underline{\omega}_{C}=p^{!}\mathcal{O}_{\mathrm{Spec}\, \mathbb{C}}$ the dualizing complex of $C$ (see \cite[Ch. V]{rd}). We thus obtain a map $$\widetilde{v}: f^{*}T_{S} \longrightarrow \underline{\omega}_{C}\otimes H^{2}(S,\mathcal{O}_{S})[-1]\simeq \omega_{C} \otimes H^{2}(S,\mathcal{O}_{S}).$$ By the properties of dualizing complexes, we have $$\underline{\omega}_{C}\otimes H^{2}(S,\mathcal{O}_{S})[-1] = \underline{\omega}_{C}\otimes p^{*}(H^{2}(S,\mathcal{O}_{S}))[-1] \simeq p^{!}(H^{2}(S,\mathcal{O}_{S})[-1]),$$ so we get a canonical morphism $$f^{*}T_{S} \longrightarrow  p^{!}(H^{2}(S,\mathcal{O}_{S})[-1])$$ which induces, by applying $\mathbb{R}p_{*}$ and composing with the adjunction map $\mathbb{R}p_{*} p^{!} \rightarrow \mathrm{Id}$, a canonical map $$\xymatrix{\widetilde{\alpha}: \mathbb{R}\Gamma (C, f^{*}T_{S}) \simeq \mathbb{R}p_{*} (f^{*}T_{S}) \ar[r]^-{\mathbb{R}p_{*}(\widetilde{v})} & \mathbb{R}p_{*}(\omega_{C} \otimes H^{2}(S,\mathcal{O}_{S}))\simeq \mathbb{R}p_{*} p^{!}(H^{2}(S,\mathcal{O}_{S})[-1]) \ar[r] &  H^{2}(S,\mathcal{O}_{S})[-1]}.$$ Since $\mathbb{R}\Gamma$ is a triangulated functor, to get a unique induced map $$\alpha: \mathbb{R}\Gamma (C, \textrm{Cone}(\mathbb{T}_{C} \rightarrow f^{*}T_{S})) \longrightarrow H^{2}(S,\mathcal{O}_{S})[-1]$$ it will be enough to observe that $\mathrm{Hom}_{\textrm{D}(\mathbb{C})}(\mathbb{R}p_{*}\mathbb{T}_{C}[1], H^{2}(S,\mathcal{O}_{S})[-1]) = 0$ (which is obvious since $\mathbb{R}p_{*}\mathbb{T}_{C}[1]$ lives in degrees $[-1,0]$, while $H^{2}(S,\mathcal{O}_{S})[-1]$ in degree $1$), and to prove the following

\begin{lem} \label{lemma}The composition $$\xymatrix{\mathbb{R}p_{*}\mathbb{T}_{C} \ar[r] & \mathbb{R}p_{*} f^{*}T_{S} \ar[r]^-{\mathbb{R}p_{*}(\widetilde{v})} & \mathbb{R}p_{*} (\omega_{C} \otimes H^{2}(S,\mathcal{O}_{S}))} $$ vanishes in the derived category $D(\mathbb{C})$.
\end{lem}

\noindent \textbf{Proof.} If $C$ is smooth, the composition $$\xymatrix{\mathbb{T}_{C} \ar[r] & f^{*}T_{S} \ar[r]^-{f^{*}\psi} & f^{*}\Omega_{S}^{1}\otimes H^{2}(S,\mathcal{O}_{S}) \ar[r]^-{s\otimes \textrm{id}} & \Omega_{C}^{1}\otimes H^{2}(S,\mathcal{O}_{S}) }$$ is obviously zero, since $\mathbb{T}_{C}\simeq T_{C}$ in this case, and a curve has no $2$-forms. For a general prestable $C$, we proceed as follows. Let's consider the composition  $$\xymatrix{\theta :  \mathbb{T}_{C} \ar[r] & f^{*}T_{S} \ar[r]^-{f^{*}\psi} & f^{*}\Omega_{S}^{1}\otimes H^{2}(S,\mathcal{O}_{S}) \ar[r]^-{s\otimes \textrm{id}} & \Omega_{C}^{1}\otimes H^{2}(S,\mathcal{O}_{S}) \ar[r]^-{t\otimes \textrm{id}} & \omega_{C} \otimes H^{2}(S,\mathcal{O}_{S}):=\mathcal{L}. }$$
On the smooth locus of $C$, $\mathcal{H}^{0}(\theta)$ is zero (by the same argument used in the case $C$ smooth), hence the image of $\mathcal{H}^{0}(\theta): \mathcal{H}^{0}(\mathbb{T}_{C})\simeq T_{C} \rightarrow \mathcal{L}$ is a torsion subsheaf of the line bundle $\mathcal{L}$. But $C$ is Cohen-Macaulay, therefore this image is $0$, i.e. $\mathcal{H}^{0}(\theta)=0$; and, obviously, $\mathcal{H}^{i}(\theta)=0$ for any $i$ (i.e. for $i=1$). Now we use the hypercohomology spectral sequences $$H^{p}(C, \mathcal{H}^{q}(\mathbb{T}_{C})) \Rightarrow \mathbb{H}^{p+q}(C, \mathbb{T}_{C}) \simeq H^{p+q}(\mathbb{R}\Gamma (C,\mathbb{T}_{C})), $$ $$ H^{p}(C, \mathcal{H}^{q}(\mathcal{L}[0])) \Rightarrow \mathbb{H}^{p+q}(C, \mathcal{L}[0]) \simeq H^{p+q}(\mathbb{R}\Gamma (C,\mathcal{L}[0])) \simeq H^{p+q}(C, \mathcal{L}),$$ to conclude that the induced maps $$H^{i}(\mathbb{R}\Gamma (\theta)): H^{i}(\mathbb{R}\Gamma (C,\mathbb{T}_{C})) \longrightarrow H^{i}(\mathbb{R}\Gamma (C,\mathcal{L}))\simeq H^{i}(C,\mathcal{L})$$ are zero for all $i$'s. Since $\mathbb{C}$ is a field, we deduce that the map $\mathbb{R}\Gamma (\theta)= \mathbb{R}p_{*}(\theta)$ is zero in $D(\mathbb{C})$ as well.  \hfill $\Box$ \\

By the Lemma above, we have therefore obtained an induced map $$\alpha: \mathbb{R}\Gamma (C, \textrm{Cone}(\mathbb{T}_{C} \rightarrow f^{*}T_{S})) \longrightarrow H^{2}(S,\mathcal{O}_{S})[-1].$$ Now, O-M-P-T reduced obstruction space is defined as $\ker \,H^{1}(\alpha)$. \\ Moreover, again by Lemma \ref{lemma}, we have an induced map $$v: \mathbb{R}p_{*}\textrm{Cone}(\mathbb{T}_{C} \rightarrow f^{*}T_{S}) \longrightarrow \mathbb{R}p_{*}(\omega_{C} \otimes H^{2}(S,\mathcal{O}_{S})),$$ and, since the canonical map $$\mathbb{R}p_{*}(\omega_{C} \otimes H^{2}(S,\mathcal{O}_{S})) \longrightarrow H^{2}(S,\mathcal{O}_{S})[-1]$$ obviously induces an isomorphism on $H^{1}$, we have that O-M-P-T reduced obstruction space is \emph{also} the kernel of the map $$H^{1}(v): H^{1}(\mathbb{R}\Gamma (C,\textrm{Cone}(\mathbb{T}_{C} \rightarrow f^{*}T_{S})) \longrightarrow H^{1}(C,\omega_{C} \otimes H^{2}(S,\mathcal{O}_{S})).$$

\noindent The following result proves the nontriviality of O-M-P-T reduced obstruction space.

\begin{prop} \label{omptnonzero}
If, as we are supposing, $\beta \neq 0$,the maps $H^1 (v)$, $H^{1}(\alpha)$, $H^{1}(\widetilde{\alpha})$, and $H^{1}(\mathbb{R}p_{*}(\widetilde{v}))$ are all nontrivial, hence surjective. 
\end{prop}

\noindent \textbf{Proof} 
The non vanishing of $H^{1}(\mathbb{R}p_{*}(\widetilde{v}))$ obviously implies all other nonvanishing statements, and the nonvanishing of $H^{1}(\mathbb{R}p_{*}(\widetilde{v}))$ is an immediate consequence of the following\footnote{We thank R. Pandharipande for pointing out this statement, of which we give here our proof.}.

\begin{lem} \label{pandha} Since the curve class $\beta \neq 0$, the map $$H^{1}(t\circ s): H^{1}(C, f^{*}\Omega_{S}^{1}) \longrightarrow H^{1}(C, \omega_{C})$$ is nonzero (hence surjective).
\end{lem} 

\noindent \textsf{Proof of Lemma.} By \cite[Cor. 2.3]{bm}, $\beta \neq 0$ implies nontriviality of the map $df: f^{*}\Omega_{S}^{1} \rightarrow \Omega_{C}^{1}$. But $S$ is a smooth surface and $C$ a prestable curve, hence in the short exact sequence $$\xymatrix{f^{*}\Omega_{S}^{1}\ar[r]^-{s} & \Omega_{C}^{1} \ar[r] & \Omega_{C/S}^{1}\rightarrow 0}$$ the sheaf of relative differentials $\Omega_{C/S}^{1}$ is concentrated at the (isolated, closed) singular points and thus its $H^{1}$ vanishes. Therefore the map $$H^{1}(s): H^{1}(C, f^{*}\Omega^{1}_{S}) \longrightarrow  H^{1}(C,\Omega_{C}^{1})$$ is surjective. The same argument yields surjectivity, hence nontriviality (since $H^{1}(C,\omega_{C})$ has dimension $1$ over $\mathbb{C}$), of the map $H^{1}(t): H^{1}(C,\Omega_{C}^{1}) \rightarrow H^{1}(C,\omega_{C})$, by observing that, on the smooth locus of $C$, $\Omega_{C}^{1}\simeq \omega_{C}$ and $H^{1}(t)$ is the induced isomorphism. In particular, $H^{1}(C,\Omega_{C}^{1}) \neq 0$.  Therefore both $H^1(s)$ and $H^{1}(t)$ are non zero and surjective, so the same is true of their composition. \,\,\,\, $\diamondsuit$ \hfill $\Box$\\

\subsection{The canonical projection $\mathbb{R}\mathbf{Pic}(S) \longrightarrow \mathbb{R}\mathrm{Spec}(\mathrm{Sym}(H^{0}(S,K_{S})[1]))$}\label{projection}

In this subsection we identify canonically the derived Picard stack $\mathbb{R}\mathbf{Pic}(S)$ of a $K3$-surface with  $\mathbf{Pic}(S) \times \mathbb{R}\mathrm{Spec}(\mathrm{Sym} (H^{0}(S, K_{S})[1]))$, where $K_{S}:=\Omega^2_{S}$ is the canonical sheaf of $S$; this allows us to define the canonical map $\textrm{pr}_{\textrm{der}}: \mathbb{R}\mathbf{Pic}(S) \longrightarrow \mathbb{R}\mathrm{Spec}(\mathrm{Sym}(H^{0}(S,K_{S})[1]))$ which is the last ingredient we will need to define the \emph{reduced} derived stack $\mathbb{R}\overline{\mathbf{M}}^{\textrm{red}}_{g}(S;\beta)$ of stable maps of genus $g$ and class $\beta$ to $S$ in the next subsection.\\

In the proof of the next Proposition, we will need the following elementary result  (which holds true for $k$ replaced by any semisimple ring, or $k$ replaced by a hereditary commutative ring and $E$ by a bounded above complex of free modules)

\begin{lem}\label{weirdmap} Let $k$ be a field and $E$ be a bounded above complex of $k$-vector spaces. Then there is a canonical map $E \rightarrow E_{<0}$ in the derived category $\mathrm{D}(k)$, such that the obvious composition $$E_{<0} \longrightarrow E \longrightarrow E_{<0}$$ is the identity.
\end{lem}

\noindent \textbf{Proof.} Any splitting $p$ of the map of $k$-vector spaces $$\ker (d_{0}: E_{-1}\rightarrow E_{0}) \hookrightarrow E_{-1}$$ yields a map $\overline{p}:E \rightarrow E_{<0}$ in the category $\textrm{Ch}(k)$ of complexes of $k$-vector spaces. To see that different splittings $p$ and $q$ gives the same map in the derived category $\mathrm{D}(k)$, we consider the canonical exact sequences of complexes $$\xymatrix{0 \rightarrow  E_{<0} \ar[r] & E \ar[r] & E_{\geq 0}  \rightarrow 0}$$ and apply $\mathrm{Ext}^{0}(-,E_{<0})$, to get an exact sequence $$\xymatrix{\mathrm{Ext}^{0}(E_{\geq 0},E_{<0}) \ar[r]^-{a} & \mathrm{Ext}^{0}(E,E_{<0}) \ar[r]^-{b} & \mathrm{Ext}^{0}(E_{< 0},E_{<0}). }$$ Now, the the class of the difference $(\overline{p}-\overline{q})$ in $\mathrm{Hom}_{\textrm{D}(k)}(E,E_{<0})= \mathrm{Ext}^{0}(E,E_{<0})$ is in the kernel of $b$, so it is enough to show that $\mathrm{Ext}^{0}(E_{\geq 0},E_{<0}) = 0$. But $E_{\geq 0}$ is a bounded above complex of projectives, therefore (e.g. \cite[Cor. 10.4.7]{we}) $\mathrm{Ext}^{0}(E_{\geq 0},E_{<0}) = 0$ is a quotient of $\mathrm{Hom}_{\mathrm{Ch}(k)}(E_{\geq 0},E_{<0})$ which obviously consists of the zero morphism alone. \hfill $\Box$ \\

\begin{prop} \label{gamma}Let $G$ be a derived group stack locally of finite presentation over  a field $k$, $\mathrm{e}: \mathrm{Spec} \, k \rightarrow G$ its identity section, and $\mathfrak{g}:=\mathbb{T}_{\textrm{e}}G$. Then there is a canonical map in $\mathrm{Ho}(\mathbf{dSt}_{k})$ $$\gamma(G): \tr (G) \times \mathbb{R}\mathrm{Spec} (A) \longrightarrow G$$ where $A:= k\oplus (\mathfrak{g}^{\vee})_{< 0}$ is the commutative differential non-positively graded $k$-algebra which is the trivial square zero extension of $k$ by the complex of $k$-vector spaces $(\mathfrak{g}^{\vee})_{< 0}$.
\end{prop}

\noindent \textbf{Proof.} First observe that $\mathbb{R}\mathrm{Spec} (A)$ has a canonical $k$-point $x_0: \mathrm{Spec} \, k \rightarrow \mathbb{R}\mathrm{Spec} (A)$, corresponding to the canonical projection $A\rightarrow k$. By definition of the derived cotangent complex of a derived stack (\cite[1.4.1]{hagII}), giving a map $\alpha$ such that $$\xymatrix{\mathbb{R}\mathrm{Spec} (A) \ar[rr]^-{\alpha} & & G \\ & \mathrm{Spec} \, k \ar[ul]^-{x_0} \ar[ur]_-{\textrm{e}} &}$$ is equivalent to giving a morphism in the derived category of complex of $k$-vector spaces $$\alpha ': \mathbb{L}_{G,\, \textrm{e}} \simeq \mathfrak{g}^{\vee} \longrightarrow (\mathfrak{g}^{\vee})_{< 0}.$$ Since $k$ is a field, we may take as $\alpha'$ the canonical map provided by Lemma \ref{weirdmap}, and define $\gamma(G)$ as the composition $$\xymatrix{\tr (G) \times \mathbb{R}\mathrm{Spec} (A) \ar[r]^-{j\times \textrm{id}} & G \times \mathbb{R}\mathrm{Spec} (A) \ar[r]^-{\textrm{id}\times \alpha'} & G\times G \ar[r]^-{\mu} & G}$$ where $\mu$ is the product in $G$.
\hfill $\Box$ \\

\begin{prop} \label{gammak3} Let $S$ be a $K3$ surface over $k=\mathbb{C}$, and $G:=\mathbb{R}\mathbf{Pic}(S)$ its derived Picard group stack. Then the map $\gamma(G)$ defined in (the proof of) Proposition \ref{gamma} is an isomorphism $$\gamma_{S}:= \gamma(\mathbb{R}\mathbf{Pic}(S)) : \xymatrix{\mathbf{Pic}(S) \times \mathbb{R}\mathrm{Spec}(\mathrm{Sym} (H^{0}(S, K_{S})[1])) \ar[r]^-{\sim} & \mathbb{R}\mathbf{Pic}(S)}$$ in $\mathrm{Ho}(\dst)$, where $K_{S}$ denotes the canonical bundle on $S$.  
\end{prop}

\noindent \textbf{Proof.} Since $G:=\mathbb{R}\mathbf{Pic}(S)$ is a derived group stack, $\gamma(G)$ is an isomorphism if and only if it induces an isomorphism on truncations, and it is \'etale at $\textrm{e}$, i.e. the induced map $$\mathbb{T}_{(\tr (\textrm{e}), x_{0})} (\gamma(G)): \mathbb{T}_{(\tr (\textrm{e}), x_{0})}(\tr (G) \times \mathbb{R}\mathrm{Spec} (A)) \longrightarrow \mathbb{T}_{\textrm{e}}(G)$$ is an isomorphism in the derived category $\mathrm{D}(k)$, where $x_0$ is the canonical canonical $k$-point $\mathrm{Spec} \, \mathbb{C} \rightarrow \mathbb{R}\mathrm{Spec} (A)$, corresponding to the canonical projection $A\rightarrow \mathbb{C}$. Since $\pi_0 (A)\simeq \mathbb{C}$,  $\tr (\gamma(G))$ is an isomorphism of stacks. So we are left to showing that $\gamma(G)$ induces an isomorphism between tangent spaces. Now, $$ \mathfrak{g}\equiv \mathbb{T}_{\textrm{e}}(G)=\mathbb{T}_{\textrm{e}}(\mathbb{R}\mathbf{Pic}(S))\simeq \mathbb{R}\Gamma(S, \mathcal{O}_{S})[1],$$ and, $S$ being a K3-surface, we have $$\mathfrak{g} \simeq \mathbb{R}\Gamma(S, \mathcal{O}_{S})[1]\simeq H^{0}(S,\mathcal{O}_{S})[1]\oplus H^{2}(S,\mathcal{O}_{S})[-1]$$ so that $$(\mathfrak{g}^{\vee})_{< 0} \simeq H^{2}(S,\mathcal{O}_{S})^{\vee}[1] \simeq H^{0}(S,K_{S})[1]$$ (where we have used Serre duality in the last isomorphism). But $H^{0}(S,K_{S})$ is free of dimension $1$, so we have a canonical isomorphism $$\mathbb{C} \oplus (\mathfrak{g}^{\vee})_{< 0} \simeq \mathbb{C} \oplus  H^{0}(S,K_{S})[1] \simeq \mathrm{Sym} (H^{0}(S, K_{S})[1])$$ in the homotopy category of commutative simplicial $\mathbb{C}$-algebras. Therefore $$\mathbb{T}_{(\tr (\textrm{e}), x_{0})}(\mathbf{Pic}(S) \times \mathbb{R}\mathrm{Spec} (A)) \simeq \mathfrak{g}_{\leq 0} \oplus \mathfrak{g}_{> 0} \simeq H^{0}(S,\mathcal{O}_{S})[1]\oplus H^{2}(S,\mathcal{O}_{S})[-1]$$ and $\mathbb{T}_{(\tr (\textrm{e}), x_{0})} (\gamma(G))$ is obviously an isomorphism (given, in the notations of the proof of Prop. \ref{gamma}, by the sum of the dual of $\alpha'$ and the canonical map $\mathfrak{g}_{\leq 0} \rightarrow \mathfrak{g}$).
 \hfill $\Box$ \\

Using Prop. \ref{gammak3}, we are now able to define the projection $\textrm{pr}_{\textrm{der}}$ of $\mathbb{R}\mathbf{Pic}(S)$ onto its full derived factor as the composite $$\xymatrix{ \mathbb{R}\mathbf{Pic}(S) \ar[r]_-{\gamma(S)^{-1}} & \mathbf{Pic}(S) \times \mathbb{R}\mathrm{Spec}(\mathrm{Sym} (H^{0}(S, K_{S})[1])) \ar[r]_-{\textrm{pr}_{2}} & \mathbb{R}\mathrm{Spec}(\mathrm{Sym} (H^{0}(S, K_{S})[1])).}$$ Note that $\textrm{pr}_{\textrm{der}}$ yields on tangent spaces the canonical projection\footnote{Recall that, if $M$ is a $\mathbb{C}$-vector space, $\mathbb{T}_{x_{0}}(\mathbb{R}\mathrm{Spec}(\textrm{Sym}(M[1])))\simeq M^{\vee}[-1]$.} $$\mathbb{T}_{\textrm{e}}(\mathbb{R}\mathbf{Pic}(S;\beta)) = \mathfrak{g} \longrightarrow \mathfrak{g}_{> 0} = \mathbb{T}_{x_{0}}(\mathbb{R}\mathrm{Spec}(\mathrm{Sym} (H^{0}(S, K_{S})[1])))\simeq H^{2}(S,\mathcal{O}_S)[-1],$$ where $x_0$ is the canonical canonical $k$-point $\mathrm{Spec} \, \mathbb{C} \rightarrow \textrm{Spec}(\mathrm{Sym} (H^{0}(S, K_{S})[1]))$, and  $$\mathfrak{g}\simeq H^{0}(S,\mathcal{O}_{S})[1]\oplus H^{2}(S,\mathcal{O}_{S})[-1].$$

\subsection{The reduced derived stack of stable maps $\mathbb{R}\overline{\mathbf{M}}^{\textrm{red}}_{g}(S;\beta)$}\label{our}

In this subsection we define the \emph{reduced} version of the derived stack of stable maps of type $(g,\beta)$ to $S$ and describe the obstruction theory it induces on its truncation $\overline{\mathbf{M}}_{g}(S;\beta)$.\\

Let us define $\delta_{1}^{\textrm{der}}(S, \beta)$ (respectively, $\delta_{1}^{(n),\,\textrm{der}}(S, \beta)$)  as the composition (see Def. \ref{delta1} and Def. \ref{delta1pointed}) $$\xymatrix{\mathbb{R}\overline{\mathbf{M}}_{g}(S;\beta) \ar@{^{(}->}[r] & \mathbb{R}\overline{\mathbf{M}}_{g}(S) \ar[r]^-{\delta_{1}(S)} & \mathbb{R}\mathbf{Pic}(S) \ar[r]^-{\textrm{pr}_{\textrm{der}}} & \mathbb{R}\mathrm{Spec}(\mathrm{Sym}(H^{0}(S,K_{S})[1]))}$$ (resp. as the composition $$\xymatrix{\mathbb{R}\overline{\mathbf{M}}_{g,n}(S;\beta) \ar@{^{(}->}[r] & \mathbb{R}\overline{\mathbf{M}}_{g,n}(S) \ar[r]^-{\delta^{(n)}_{1}(S)} & \mathbb{R}\mathbf{Pic}(S) \ar[r]^-{\textrm{pr}_{\textrm{der}}} & \mathbb{R}\mathrm{Spec}(\mathrm{Sym}(H^{0}(S,K_{S})[1]))} \, \,\, ).$$ 

\begin{df} \label{dereduced}
\begin{itemize} 
\item The \emph{reduced} derived stack of stable maps of genus $g$ and class $\beta$ to $S$ $\mathbb{R}\overline{\mathbf{M}}^{\textrm{red}}_{g}(S;\beta)$ is defined by the following homotopy-cartesian square in $\dst$
$$\xymatrix{ \mathbb{R}\overline{\mathbf{M}}^{\textrm{red}}_{g}(S;\beta) \ar[d] \ar[r] & \mathbb{R}\overline{\mathbf{M}}_{g}(S;\beta) \ar[d]^-{\delta_{1}^{\textrm{der}}(S,\beta)} \\
\mathrm{Spec}\, \mathbb{C} \ar[r] & \mathbb{R}\mathrm{Spec}(\mathrm{Sym}(H^{0}(S,K_{S})[1])) }$$
\item The \emph{reduced} derived stack of $n$-pointed stable maps of genus $g$ and class $\beta$ to $S$ $\mathbb{R}\overline{\mathbf{M}}^{\textrm{red}}_{g,n}(S;\beta)$ is defined by the following homotopy-cartesian square in $\dst$
$$\xymatrix{ \mathbb{R}\overline{\mathbf{M}}^{\textrm{red}}_{g,n}(S;\beta) \ar[d] \ar[r] & \mathbb{R}\overline{\mathbf{M}}_{g,n}(S;\beta) \ar[d]^-{\delta_{1}^{(n), \textrm{der}}(S,\beta)} \\
\mathrm{Spec}\, \mathbb{C} \ar[r] & \mathbb{R}\mathrm{Spec}(\mathrm{Sym}(H^{0}(S,K_{S})[1])) }$$
\end{itemize}
\end{df}

Since the truncation functor $\tr$ commutes with homotopy fiber products and $$\tr (\mathbb{R}\mathrm{Spec}(\mathrm{Sym}(H^{0}(S,K_{S})[1])))\simeq \mathrm{Spec} \, \mathbb{C},$$ we get $$\tr (\mathbb{R}\overline{\mathbf{M}}^{\textrm{red}}_{g}(S;\beta))\simeq \overline{\mathbf{M}}_{g}(S;\beta)$$ i.e. $\mathbb{R}\overline{\mathbf{M}}^{\textrm{red}}_{g}(S;\beta)$ is a \emph{derived extension} (Def. \ref{derext}) of the usual stack of stable maps 
of type $(g, \beta)$ to $S$, different from $\mathbb{R}\overline{\mathbf{M}}_{g}(S;\beta)$. Similarly in the pointed case.\\

We are now able to compute the obstruction theory induced, according to \S \ref{derobs}, by the closed immersion $j_{\mathrm{red}}: \overline{\mathbf{M}}_{g}(S;\beta) \hookrightarrow \mathbb{R}\overline{\mathbf{M}}^{\textrm{red}}_{g}(S;\beta)$. We leave to the reader the straightforward modifications for the pointed case.\\  

By applying Proposition \ref{inducedobstruction} to the derived extension $ \mathbb{R}\overline{\mathbf{M}}^{\textrm{red}}_{g}(S;\beta)$ of $\overline{\mathbf{M}}_{g}(S;\beta)$, we get an obstruction theory $$j_{\mathrm{red}}^{*}\mathbb{L}_{\mathbb{R}\overline{\mathbf{M}}^{\textrm{red}}_{g}(S;\beta)} \longrightarrow \mathbb{L}_{\overline{\mathbf{M}}_{g}(S;\beta)}$$ that we are now going to describe.\\ Let $$\rho:\mathbb{R}\overline{\mathbf{M}}^{\textrm{red}}_{g}(S;\beta) \longrightarrow \mathbb{R}\overline{\mathbf{M}}_{g}(S;\beta)$$ be the canonical map. Since $\mathbb{R}\overline{\mathbf{M}}^{\textrm{red}}_{g}(S;\beta)$ is defined by the homotopy pullback diagram in Def. \ref{dereduced}, we get an isomorphism in the derived category of $\mathbb{R}\overline{\mathbf{M}}^{\textrm{red}}_{g}(S;\beta)$ $$ \rho^{*}(\mathbb{L}_{\mathbb{R}\overline{\mathbf{M}}_{g}(S;\beta)/\mathbb{R}\mathrm{Spec}(\mathrm{Sym}(H^{0}(S,K_{S})[1]))}) \simeq \mathbb{L}_{\mathbb{R}\overline{\mathbf{M}}^{\textrm{red}}_{g}(S;\beta)}.$$ We will show below that  $\mathbb{R}\overline{\mathbf{M}}^{\textrm{red}}_{g}(S;\beta)$ is \emph{quasi-smooth} so that, by Corollary \ref{cor}, $$j_{\mathrm{red}}^{*}\mathbb{L}_{\mathbb{R}\overline{\mathbf{M}}^{\textrm{red}}_{g}(S;\beta)} \longrightarrow \mathbb{L}_{\overline{\mathbf{M}}_{g}(S;\beta)}$$ is indeed a \emph{perfect} obstruction theory on $\overline{\mathbf{M}}_{g}(S;\beta)$. Now, for any $\mathbb{C}$-point $\mathrm{Spec}\, \mathbb{C} \rightarrow \mathbb{R}\overline{\mathbf{M}}_{g}(S;\beta)$, corresponding to a stable map $(f:C\rightarrow S)$ of type $(g,\beta)$, we get a distinguished triangle $$\mathbb{L}_{\mathbb{R}\mathrm{Spec}(\mathrm{Sym}(H^{0}(S,K_{S})[1])),\, x_{0}} \longrightarrow \mathbb{L}_{\mathbb{R}\overline{\mathbf{M}}_{g}(S;\beta) ,\,(f:C\rightarrow S)} \longrightarrow \mathbb{L}_{\mathbb{R}\overline{\mathbf{M}}^{\textrm{red}}_{g}(S;\beta), \, (f:C\rightarrow S)}$$ (where we have denoted by $(f:C\rightarrow S)$ also the induced $\mathbb{C}$-point of $\mathbb{R}\overline{\mathbf{M}}^{\textrm{red}}_{g}(S;\beta)$: recall that a derived stack and its truncation have the same classical points, i.e. points with values in usual commutative $\mathbb{C}$-algebras) in the derived category of complexes of $\mathbb{C}$-vector spaces. By dualizing, we get that the tangent complex $$\mathbb{T}^{\textrm{red}}_{(f:C\rightarrow S)}:=\mathbb{T}_{(f:C\rightarrow S)}(\mathbb{R}\overline{\mathbf{M}}^{\textrm{red}}_{g}(S;\beta))$$ of $\mathbb{R}\overline{\mathbf{M}}^{\textrm{red}}_{g}(S;\beta)$ at the $\mathbb{C}$-point $(f:C\rightarrow S)$ of type $(g,\beta)$, sits into a distinguished triangle $$\xymatrix{\mathbb{T}^{\textrm{red}}_{(f:C\rightarrow S)} \ar[r] & \mathbb{R}\Gamma(C, \textrm{Cone}(\mathbb{T}_{C}\rightarrow f^{*}T_{S})) \ar[r]^-{\Theta_{f}} & \mathbb{R}\Gamma (S, \mathcal{O}_{S})[1] \ar[r]^-{\textrm{pr}} & H^{2}(S,\mathcal{O}_{S})[-1]},$$  where $\Theta_{f}$ is the composite $$\xymatrix{\Theta_{f} : \mathbb{R}\Gamma(C,\mathrm{Cone}(\mathbb{T}_{C}\rightarrow f^{*}T_{S})) \ar[r]^-{\mathbb{T}_{x}\mathrm{A}_{X}} & \mathbb{R}\mathrm{Hom}_{S}(\mathbb{R}f_{*}\mathcal{O}_{C}, \mathbb{R}f_{*}\mathcal{O}_{C})[1] \ar[r]^-{\mathrm{tr}_{S}} & \mathbb{R}\Gamma (S, \mathcal{O}_{S})[1],}$$ and $\textrm{pr}$ denotes the tangent map of $\textrm{pr}_{\textrm{der}}$ taken at the point $\delta_{1}(S)(f:C\rightarrow S)$. Note that the map pr obviously induces an isomorphism on $H^{1}$.\\

\subsection{Quasi-smoothness of $\mathbb{R}\overline{\mathbf{M}}^{\textrm{red}}_{g}(S;\beta)$ and comparison with O-M-P-T reduced obstruction theory.} In the case $\beta \neq 0$ is a curve class in $H^{2}(S,\mathbb{Z})$, we will prove \emph{quasi-smoothness} of the derived stack $\mathbb{R}\overline{\mathbf{M}}^{\textrm{red}}_{g}(S;\beta)$, and compare the induced obstruction theory with that of Okounkov-Maulik-Pandharipande-Thomas (see \S \ref{omptspaces} or \cite[\S 2.2]{mp} and \cite{pak3}).\\

\begin{thm} \label{quasismooth} Let $\beta \neq 0$ be a curve class in $H^{2}(S,\mathbb{Z})\simeq H_{2}(S,\mathbb{Z})$, $f:C\rightarrow S$ a stable map of type $(g,\beta)$, and $$\xymatrix{\mathbb{T}^{\textrm{red}}_{(f:C\rightarrow S)}:=\mathbb{T}_{(f:C\rightarrow S)}(\mathbb{R}\overline{\mathbf{M}}^{\textrm{red}}_{g}(S;\beta)) \ar[r] & \mathbb{R}\Gamma(C, \textrm{Cone}(\mathbb{T}_{C}\rightarrow f^{*}T_{S})) \ar[r] & H^{2}(S,\mathcal{O}_{S})[-1] }$$ the corresponding distinguished triangle. Then, 
\begin{enumerate}
\item the rightmost arrow in the triangle above induces on $H^1$ a map $$H^1(\Theta_{f}): H^{1}(C, \textrm{Cone}(\mathbb{T}_{C}\rightarrow f^{*}T_{S})) \longrightarrow H^{2}(S, \mathcal{O}_{S})$$ which is nonzero (hence surjective, since $H^{2}(S, \mathcal{O}_{S})$ has dimension $1$ over $\mathbb{C}$). Therefore the derived stack $\mathbb{R}\overline{\mathbf{M}}^{\textrm{red}}_{g}(S;\beta)$ is everywhere quasi-smooth;
\item $H^{0}(\mathbb{T}^{\textrm{red}}_{(f:C\rightarrow S)})$ (resp.  $H^{1}(\mathbb{T}^{\textrm{red}}_{(f:C\rightarrow S)})$) coincides with the reduced deformation space (resp. the reduced obstruction space) of O-M-P-T. 
\end{enumerate}
\end{thm}

\noindent \textbf{Proof.} \\
\noindent \textsf{First Proof of quasi-smoothness --} Let us prove quasi-smoothness first. It is clearly enough to prove that the composite $$\xymatrix{H^{1}(C,f^{*}T_{S}) \ar[r] & H^{1}(C,\mathrm{Cone}(\mathbb{T}_{C}\rightarrow f^{*}T_{S})) \ar[rr]^-{H^{1}(\mathbb{T}_{x}\mathrm{A}_{X})} & & \mathrm{Ext}^{2}_{S}(\mathbb{R}f_{*}\mathcal{O}_{C}, \mathbb{R}f_{*}\mathcal{O}_{C}) \ar[r]^-{H^{2}(\mathrm{tr}_{S})} & H^{2}(S, \mathcal{O}_{S})}$$ is non zero (hence surjective). Recall that $p: C\rightarrow \mathrm{Spec}\, \mathbb{C}$ and $q:S \rightarrow \mathrm{Spec}\, \mathbb{C}$ denote the structural morphisms, so that $p=q\circ f$. Now, the map $$\mathbb{R}q_{*}T_{S} \longrightarrow \mathbb{R}q_{*}\mathbb{R}f_{*}f^{*}T_{S}$$ induces a map $H^{1}(S,T_{S})\rightarrow H^{1}(C,f^{*}T_{S})$, and by Proposition \ref{finallyatiyah} and Remark \ref{atiyah2}, the following diagram commutes $$\xymatrix{H^{1}(S,T_{S}) \ar[d] \ar[rr]^-{<-, \textrm{at}_{\mathbb{R}f_{*}\mathcal{O}_{C}}>} &  & \mathrm{Ext}^{2}_{S}(\mathbb{R}f_{*}\mathcal{O}_{C}, \mathbb{R}f_{*}\mathcal{O}_{C}) \\ H^{1}(C,f^{*}T_{S}) \ar[rr] & & H^{1}(C,\textrm{Cone}(\mathbb{T}_{C} \rightarrow f^{*}T_{S})) \ar[u]_-{H^{1}(\mathbb{T}_{x}\mathrm{A}_{X})} }.$$ So, we are reduced to proving that the composition 
$$\xymatrix{a: H^{1}(S,T_{S}) \ar[rr]^-{<-, \textrm{at}_{\mathbb{R}f_{*}\mathcal{O}_{C}}>} &  & \mathrm{Ext}^{2}_{S}(\mathbb{R}f_{*}\mathcal{O}_{C}, \mathbb{R}f_{*}\mathcal{O}_{C}) \ar[r]^-{H^2(\textrm{tr})} & H^2(S,\mathcal{O}_{S})}$$ does not vanish. But, since the first Chern class is the trace of the Atiyah class, this composition acts as follows (on maps in the derived category of $S$) $$\xymatrix{(\xi : \mathcal{O}_{S} \rightarrow T_{S}[1]) \ar[rr] & & (a(\xi): \mathcal{O}_{S} \ar[r]^-{c_{1}\otimes \xi} & \Omega^{1}_{S}\otimes T_{S}[2] \ar[r]^-{<-,->} & \mathcal{O}_{S}[2]) }$$ where $$c_{1}:= c_{1}(\mathbb{R}f_{*}\mathcal{O}_{C}):\mathcal{O}_{S} \longrightarrow \Omega_{S}^{1}[1]$$ is the first Chern class of the perfect complex $\mathbb{R}f_{*}\mathcal{O}_{C}$. What we have said so far, is true for any smooth complex projective scheme $X$ in place of $S$. We now use the fact that $S$ is a $K3$-surface. Choose a non zero section $\sigma: \mathcal{O}_{S} \rightarrow \Omega_{S}^{2}$ of the canonical bundle, and denote by $\xymatrix{\varphi_{\sigma}: \Omega^1_S \ar[r]^-{\sim} & T_{S}}$ the induced isomorphism. A straightforward linear algebra computation shows then that the composition $$\xymatrix{\mathcal{O}_{S} \ar[rrr]^-{((\varphi_{\sigma}\circ c_1)\wedge \xi )\otimes \sigma}  &&& (T_{S}\wedge T_{S} \otimes \Omega_{S}^{2})[2] \ar[rr]^-{<-,->[2]} && \mathcal{O}_{S}[2] }$$ coincides with $a(\xi)$. But, since $\beta \neq 0$, we have that $c_1\neq 0$. $\sigma$ is nondegenerate, so this composition cannot vanish for all $\xi$, and we conclude.\\

\noindent \textsf{Second Proof of quasi-smoothness --} Let us give an alternative proof of quasi-smoothness. By Serre duality, passing to dual vector spaces and maps, we are left to proving that the composite $$\xymatrix{H^0(S,\Omega^{2}_{S}) \ar[r]^-{\textrm{tr}^{\vee}} & \mathrm{Ext}_{S}^{0}(\mathbb{R}\underline{\mathrm{Hom}}(\mathbb{R}f_{*}\mathcal{O}_{C}, \mathbb{R}f_{*}\mathcal{O}_{C}), \Omega_{S}^{2}) \ar[r]^-{\tau^{\vee}} & \mathrm{Ext}^{0}(\mathbb{R}f_{*}f^{*}T_{S}[-1], \Omega_{S}^{2}) }$$ is non zero. So it is enough to prove that the map obtained by further composing to the left with the adjunction map $$\mathrm{Ext}^{0}(\mathbb{R}f_{*}f^{*}T_{S}[-1], \Omega_{S}^{2}) \longrightarrow \mathrm{Ext}^{0}(T_{S}[-1], \Omega_{S}^{2})$$ is nonzero. But this new composition acts as follows $$H^0(S,\Omega^{2}_{S}) \ni (\sigma:\mathcal{O}_{S} \rightarrow \Omega_{S}^{2}) \mapsto (\sigma\circ \mathrm{tr}) \mapsto (\sigma \circ \mathrm{tr} \circ \mathrm{at})=(\sigma \circ c_{1}(\mathbb{R}f_{*}\mathcal{O}_{C})) \in \mathrm{Ext}^{0}(T_{S}[-1], \Omega_{S}^{2})$$ where $\mathrm{at}: T_{S}[-1] \rightarrow \mathbb{R}\underline{\mathrm{Hom}}(\mathbb{R}f_{*}\mathcal{O}_{C}, \mathbb{R}f_{*}\mathcal{O}_{C})$ is the Atiyah class of $\mathbb{R}f_{*}\mathcal{O}_{C}$ (see  Proposition \ref{finallyatiyah} and Remark \ref{atiyah2}). Since $\beta \neq 0$, we have $c_{1}(\mathbb{R}f_{*}\mathcal{O}_{C})\neq 0$, and we conclude.\\

\noindent\textsf{Proof of the comparison --} Let us move now to the second point of Thm. \ref{quasismooth}, i.e. the comparison statement about deformations and obstructions spaces. First of all it is clear that, for any $\beta$, $$H^{0}(\mathbb{T}_{\textrm{red},\, (f:C\rightarrow S)}) \simeq H^{0}(\mathbb{T}_{\mathbb{R}\overline{\mathbf{M}}_{g}(S;\beta),\,(f:C\rightarrow S)}) \simeq H^{0}(C, \textrm{Cone}(\mathbb{T}_{C}\rightarrow f^{*}T_{S}))$$ therefore our deformation space is the same as O-M-P-T's one. Let us then concentrate on obstruction spaces. \\

We begin by noticing the following fact

\begin{lem}\emph{There is a canonical morphism in $\mathrm{D}(\mathbb{C})$ $$\nu: \mathbb{R}p_{*}\omega_{C}\otimes^{\mathbb{L}}H^{2}(S,\mathcal{O}_{S}) \longrightarrow \mathbb{R}q_{*}\mathcal{O}_{S}[1]$$ inducing an isomorphism on $H^{1}$.}
\end{lem}

\noindent\textsf{Proof of Lemma.} To ease notation we will simply write $\otimes$ for $\otimes^{\mathbb{L}}$. Recall that $p: C\rightarrow \mathrm{Spec}\, \mathbb{C}$ and $q:S \rightarrow \mathrm{Spec}\, \mathbb{C}$ denote the structural morphisms, so that $p=q\circ f$. Since $S$ is a $K3$-surface, the canonical map $$\mathcal{O}_{S}\otimes H^{0}(S,\Omega^{2}_{S})\longrightarrow  \Omega_{S}^{2}$$ is an isomorphism. Since $f^{!}$ preserves dualizing complexes, $\underline{\omega}_{S}\simeq \Omega_{S}^{2}[2]$ and $\underline{\omega}_{C}\simeq \omega_{C}[1]$, we have $$\omega_{C}\simeq f^{!}\Omega_{S}^{2}[1] \simeq f^{!}(\mathcal{O}_{S}\otimes H^{0}(S,\Omega^{2}_{S}))[1].$$ By applying $\mathbb{R}p_{*}$ and using the adjunction map $\mathbb{R}f_{*}f^{!} \rightarrow \mathrm{Id}$, we get a map
$$\mathbb{R}p_{*}\omega_{C} \simeq \mathbb{R}q_{*}\mathbb{R}f_{*}\omega_{C}\simeq \mathbb{R}q_{*}\mathbb{R}f_{*}f^{!}(\mathcal{O}_{S}[1]\otimes H^{0}(S,\Omega^{2}_{S})) \rightarrow \mathbb{R}q_{*}(\mathcal{O}_{S}[1]\otimes H^{0}(S,\Omega^{2}_{S})) \simeq \mathbb{R}q_{*}\mathcal{O}_{S}[1]\otimes H^{0}(S,\Omega^{2}_{S})$$ (the last isomorphism being given by projection formula). Tensoring this map by $H^{0}(S,\Omega^{2}_{S})^{\vee}$ $\simeq H^{2}(S,\mathcal{O}_{S})$ (a canonical isomorphism by Serre duality), and using the canonical evaluation map $V\otimes V^{\vee} \rightarrow \mathbb{C}$ for a $\mathbb{C}$-vector space $V$, we get the desired canonical map $$\nu: \mathbb{R}p_{*}\omega_{C}\otimes H^{2}(S,\mathcal{O}_{S}) \longrightarrow \mathbb{R}q_{*}\mathcal{O}_{S}[1].$$ The isomorphism on $H^{1}$ is obvious since the trace map $\mathrm{R}^{1}p_{*}\omega_{C} \rightarrow \mathbb{C}$ is an isomorphism ($C$ is geometrically connected). \,\,\, $\diamondsuit$ \\

If $\xymatrix{\sigma: \mathcal{O}_{S} \ar[r]^-{\sim} & \Omega_{S}^{2}}$ is a nonzero element in $H^{0}(S,\Omega_S^2)$, and $\varphi_{\sigma}: T_{S} \simeq \Omega_{S}^{1}$ the induced isomorphism, the previous Lemma gives us  an induced map $$\nu(\sigma): \mathbb{R}p_{*}\omega_{C} \longrightarrow \mathbb{R}q_{*}\mathcal{O}_{S}[1],$$ and an induced isomorphism $$\xymatrix {H^1(\nu(\sigma))=:\nu_{\sigma}:H^{1}(C,\omega_{C}) \ar[r]^-{\sim} & H^{2}(S,\mathcal{O}_{S}).}$$ Using the same notations as in  \S \ref{omptspaces}, to prove that our reduced obstruction space $$\ker (H^1(\Theta_{f}): H^{1}(C, \textrm{Cone}(\mathbb{T}_{C}\rightarrow f^{*}T_{S})) \longrightarrow H^{2}(S, \mathcal{O}_{S}))$$ coincides with O-M-P-T's one, it will be enough to show that the following diagram is commutative $$\xymatrix{H^{1}(C, f^{*}T_{S}) \ar[d]_-{\textrm{can}} \ar[r]^-{\textrm{can}} & H^{1}(\mathbb{R}\Gamma (C,\textrm{Cone}(\mathbb{T}_{C} \rightarrow f^{*}T_{S}))) \ar[dd]^-{H^1(\Theta_{f})} \\  H^{1}(\mathbb{R}\Gamma (C,\textrm{Cone}(\mathbb{T}_{C} \rightarrow f^{*}T_{S}))) \ar[d]_-{H^{1}(v)} & \\  H^{1}(C,\omega_{C}) \ar[r]^-{\sim}_-{\nu_{\sigma}} & H^{2}(S,\mathcal{O}_{S}).}$$ But this follows from the commutativity of 
$$\xymatrix{\mathbb{R}p_{*} f^{*}T_{S}[-1] \ar[d]_-{\textrm{id}} \ar[r]^-{\mathbb{R}p_{*}(\varphi_{\sigma})} & \mathbb{R}p_{*}f^{*}\Omega_{S}^{1}[-1] \ar[r]^-{\mathbb{R}p_{*}(s)} &  \mathbb{R}p_{*}\Omega_{C}^{1}[-1] \ar[r]^-{\mathbb{R}p_{*}(t)} & \mathbb{R}p_{*}\omega_{C}[-1] \ar[d]^-{\nu(\sigma)[-1]} \\  \mathbb{R}p_{*}f^{*}T_{S}[-1] \ar[r]_-{\mathbb{T}_{x}\mathrm{A}_{X}} & \mathbb{R}q_{*}\mathbb{R}\underline{\mathrm{Hom}}_{S}(\mathbb{R}f_{*}\mathcal{O}_{C}, \mathbb{R}f_{*}\mathcal{O}_{C}) \ar[rr]_-{\textrm{tr}} & & \mathbb{R}q_{*}\mathcal{O}_{S} }$$ whose verification is left to the reader.  \\

\hfill $\Box$ \\

\begin{rmk}\emph{Note that by Lemma \ref{omptnonzero}, the second assertion of Theorem \ref{quasismooth} implies the first one. Nonetheless, we have preferred to give an independent proof of the quasismoothness of $\mathbb{R}\overline{\mathbf{M}}^{\textrm{red}}_{g}(S;\beta)$ because we find it  conceptually more relevant than the comparison with O-M-P-T, meaning that quasi-smoothness alone would in any case imply the existence of  \emph{some} perfect reduced obstruction theory on $\overline{\mathbf{M}}_{g}(S;\beta)$, regardless of its comparison with the one introduced and studied by O-M-P-T. \\ Moreover, we could only find in the literature a definition of O-M-P-T \emph{global} reduced obstruction theory (relative to $\mathbf{M}^{\textrm{pre}}_{g}$) with values in the $\tau_{\geq -1}$-truncation of the cotangent complex of the stack of stable maps\footnote{The reason being that the authors use factorization through the cone, and therefore the resulting obstruction theory is only well-defined, without further arguments, if one considers it as having values in such a truncation.}, that uses a result on the semiregularity map whose proof is not completely convincing (\cite[2.2, formula (14)]{mp}); on the other hand there is a clean and complete description of the corresponding pointwise tangent and obstruction spaces. Therefore, our comparison is necessarily limited to these spaces. And our construction might also be seen as establishing such a reduced \emph{global} obstruction theory - in the usual sense, i.e. with values in the full cotangent complex, and completely independent from any result on semiregularity maps.}
\end{rmk}

\bigskip

Theorem \ref{quasismooth} shows that the distinguished triangle $$\mathbb{T}^{\textrm{red}}_{(f:C\rightarrow S)}:=\mathbb{T}_{(f:C\rightarrow S)}(\mathbb{R}\overline{\mathbf{M}}^{\textrm{red}}_{g}(S;\beta)) \longrightarrow \mathbb{R}\Gamma(C, \textrm{Cone}(\mathbb{T}_{C}\rightarrow f^{*}T_{S})) \longrightarrow H^{2}(S,\mathcal{O}_{S})[-1]$$ induces isomorphisms $$H^{i}(\mathbb{T}^{\textrm{red}}_{(f:C\rightarrow S)}) \simeq H^{i}(C, \textrm{Cone}(\mathbb{T}_{C}\rightarrow f^{*}T_{S})),$$ for any $i\neq 1$, while in degree $1$, it yields a short exact sequence  $$\xymatrix{0 \rightarrow H^{1}(\mathbb{T}^{\textrm{red}}_{(f:C\rightarrow S)}) \ar[r] & H^{1}(C, \textrm{Cone}(\mathbb{T}_{C}\rightarrow f^{*}T_{S})) \ar[r] & H^{2}(S, \mathcal{O}_{S}) \rightarrow 0.} $$
So, the tangent complexes of $\mathbb{R}\overline{\mathbf{M}}^{\textrm{red}}_{g}(S;\beta)$ and $\mathbb{R}\overline{\mathbf{M}}_{g}(S;\beta)$ (hence our induced \emph{reduced} and the \emph{standard} obstruction theories) \emph{only differ at the level of} $H^{1}$ where the former is the kernel of a $1$-dimensional quotient of the latter: this is indeed the distinguished feature of a reduced obstruction theory.\\

\noindent \textbf{The pointed case -} In the pointed case, a completely analogous proof as that of Theorem \ref{quasismooth} (1), yields 

\begin{thm} \label{quasismoothpointed} Let $\beta \neq 0$ be a curve class in $H^{2}(S,\mathbb{Z})\simeq H_{2}(S,\mathbb{Z})$. The derived stack $\mathbb{R}\overline{\mathbf{M}}^{\textrm{red}}_{g,n}(S;\beta)$ of $n$-pointed stable maps of type $(g,\beta)$ is everywhere quasi-smooth, and therefore the canonical map $$j^{*} (\mathbb{L}_{\mathbb{R}\overline{\mathbf{M}}_{g,n}(X;\beta)}) \longrightarrow \mathbb{L}_{\overline{\mathbf{M}}_{g,n}(X;\beta)}$$ is a $[-1,0]$ perfect obstruction theory on $\overline{\mathbf{M}}_{g,n}(X;\beta)$.
\end{thm}

\section{Moduli of perfect complexes}

In this Section we will define and study derived versions of various stacks of perfect complexes on a smooth projective variety $X$. If $X$ is a $K3$-surface, by using the determinant map and the structure of $\mathbb{R}\mathbf{Pic}(X)$, we deduce that the derived stack of simple perfect complexes on $X$ is smooth. This result was proved with different methods by Inaba in \cite{ina}. \\
When $X$ is a Calabi-Yau $3$-fold we prove that the derived stack of simple perfect complexes (with fixed determinant) is quasi-smooth, and then use an elaboration of the map $\mathrm{A}^{(n)}_{X}: \mathbb{R}\overline{\mathbf{M}}_{g,n}(X) \longrightarrow \mathbb{R}\mathbf{Perf}(X)$  to compare the obstruction theories induced on the truncation stacks. This might be seen as a derived geometry approach to a baby, open version of the Gromov-Witten/Donaldson-Thomas comparison.\\


\begin{df}\label{simples} Let $X$ be a smooth complex projective variety, $\mathcal{L}$ a line bundle on $X$, and $x_{\mathcal{L}}: \mathrm{Spec}\, \mathbb{C} \rightarrow \mathbb{R}\mathbf{Pic}(X)$ the corresponding point.  
\begin{itemize}
\item The \emph{derived stack $\mathbb{R}\mathbf{Perf}(X)_{\mathcal{L}}$ of perfect complexes  on $X$ with fixed determinant $\mathcal{L}$} is defined by the following homotopy cartesian diagram in $\dst$ $$\xymatrix{\mathbb{R}\mathbf{Perf}(X)_{\mathcal{L}} \ar[r] \ar[d] & \mathbb{R}\mathbf{Perf}(X) \ar[d]^-{\mathrm{det}} \\ \mathrm{Spec}\, \mathbb{C} \ar[r]_{x_{\mathcal{L}}} & \mathbb{R}\mathbf{Pic}(X)}$$  We will write $\mathbb{R}\mathbf{Perf}(X)_{0}$ for $\mathbb{R}\mathbf{Perf}(X)_{\mathcal{O}_{X}}$, the \emph{derived stack of perfect complexes on $X$ with trivial determinant}. 
\item If $\mathbf{Perf}(X)^{\geq 0}$ denotes the open substack of $\mathbf{Perf}(X)$ consisting of perfect complexes $F$ on $X$ such that $\mathrm{Ext}^i(F,F)=0$ for $i < 0$, we define $\mathbb{R}\mathbf{Perf}(X)^{\geq 0}:= \phi_{\mathbb{R}\mathbf{Perf}(X)}(\mathbf{Perf}(X)^{\geq 0})$ (as a derived open substack of  $\mathbb{R}\mathbf{Perf}(X)$, see Prop. \ref{opens}). 
\item If $\mathbf{Perf}(X)^{\textrm{si}, > 0}$ denotes the open substack of $\mathbf{Perf}(X)$ consisting of perfect complexes $F$ on $X$ for which $\mathrm{Ext}^i(F,F)=0$ for $i < 0$, and the trace map $\mathrm{Ext}^0(F,F)\rightarrow H^{0}(X,\mathcal{O}_{X})\simeq \mathbb{C}$ is an isomorphism, we define $\mathbb{R}\mathbf{Perf}(X)^{\textrm{si}, > 0}:= \phi_{\mathbb{R}\mathbf{Perf}(X)}(\mathbf{Perf}(X)^{\textrm{si}, > 0})$ (as a derived open substack of  $\mathbb{R}\mathbf{Perf}(X)$, see Prop. \ref{opens}). 
\item The derived stack $\mathbb{R}\mathbf{Perf}(X)^{\geq 0}_{\mathcal{L}}$ is defined by the following homotopy cartesian diagram in $\dst$ $$\xymatrix{\mathbb{R}\mathbf{Perf}(X)^{\geq 0}_{\mathcal{L}} \ar[r] \ar[d] & \mathbb{R}\mathbf{Perf}(X)^{\geq 0} \ar[d]^-{\mathrm{det}} \\ \mathrm{Spec}\, \mathbb{C} \ar[r]_{x_{\mathcal{L}}} & \mathbb{R}\mathbf{Pic}(X)}.$$ As above, we will write $\mathbb{R}\mathbf{Perf}(X)^{\geq 0}_{0}$ for $\mathbb{R}\mathbf{Perf}(X)^{\geq 0}_{\mathcal{O}_{X}}$.
\item The derived stack $\mathcal{M}_X\equiv\mathbb{R}\mathbf{Perf}(X)^{\textrm{si}, > 0}_{\mathcal{L}}$ is defined by the following homotopy cartesian diagram in $\dst$ $$\xymatrix{\mathbb{R}\mathbf{Perf}(X)^{\textrm{si}, > 0}_{\mathcal{L}} \ar[r] \ar[d] & \mathbb{R}\mathbf{Perf}(X)^{\textrm{si}, > 0} \ar[d]^-{\mathrm{det}} \\ \mathrm{Spec}\, \mathbb{C} \ar[r]_{x_{\mathcal{L}}} & \mathbb{R}\mathbf{Pic}(X)}.$$ We will write $\mathbb{R}\mathbf{Perf}(X)^{\textrm{si}, > 0}_{0}$ for $\mathbb{R}\mathbf{Perf}(X)^{\textrm{si}, > 0}_{\mathcal{O}_{X}}$.
\end{itemize}
\end{df}

\begin{prop}\label{almostDT} Let $E$ be a perfect complex on $X$ with determinant $\mathcal{L}$, and $x_{E}:\mathrm{Spec}\, \mathbb{C} \rightarrow \mathbb{R}\mathbf{Perf}(X)_{\mathcal{L}}$ the corresponding point. The tangent complex of $\mathbb{R}\mathbf{Perf}(X)_{\mathcal{L}}$ at $x_{E}$ is $\mathbb{R}\mathrm{End}(E)_0[1]$, the shifted traceless derived endomorphisms complex of $E$ \emph{(\cite[Def. 10.1.4]{hl})} (so that $H^{i}(\mathbb{R}\mathrm{End}(E)_0) = \ker (\mathrm{tr}: \mathrm{Ext}^i(E,E) \rightarrow H^{i}(X,\mathcal{O}_X))$, for any $i$).
\end{prop}

\noindent \textbf{Proof.} Let $\mathbb{T}$ denote the tangent complex of $\mathbb{R}\mathbf{Perf}(X)_{\mathcal{L}}$ at $x_{E}$. By definition of $\mathbb{R}\mathbf{Perf}(X)_{\mathcal{L}}$, we have an exact triangle in the derived category $D(\mathbb{C})$ of $\mathbb{C}$-vector spaces 
$$\xymatrix{\mathbb{T} \ar[r] & \mathbb{R}\mathrm{End}(E)[1] \ar[r]^-{\mathrm{tr}} & \mathbb{R}\Gamma(X, \mathcal{O}_X)[1]}.$$ By using the canonical map $\mathbb{R}\Gamma(X, \mathcal{O}_X) \rightarrow \mathbb{R}\mathrm{End}(E)$ we can split this triangle, and we conclude. 
\hfill $\Box$ \\

\begin{rmk} \emph{Since $\mathbb{R}\mathbf{Perf}(X)^{\geq 0}_{\mathcal{L}}$ and $\mathbb{R}\mathbf{Perf}(X)^{\textrm{si}, > 0}_{\mathcal{L}}$ are derived open substacks of $\mathbb{R}\mathbf{Perf}(X)_{\mathcal{L}}$, Proposition \ref{almostDT} holds for their tangent complexes too.}
\end{rmk}

\subsection{On $K3$ surfaces} 
By using the derived determinant map and the derived stack of perfect complexes, we are able to give another proof of a result by Inaba (\cite[Thm. 3.2]{ina}) that generalizes an earlier work by Mukai (\cite{muk}). For simplicity, we prove this result for $K3$ surfaces, the result for a general Calabi-Yau surface being similar.\\

Let $S$ be a a smooth projective $K3$ surface, and let  $\mathbb{R}\mathbf{Perf}(S)^{\textrm{si}, > 0}$ (Def. \ref{simples}) be the open derived substack of $\mathbb{R}\mathbf{Perf}(S)$ consisting of perfect complexes $F$ on $S$ for which $\mathrm{Ext}_{S}^i(F,F)=0$ for $i < 0$, and the trace map $\mathrm{Ext}_{S}^0(F,F)\rightarrow H^{0}(S,\mathcal{O}_{S})\simeq \mathbb{C}$ is an isomorphism. The truncation $\mathbf{Perf}(S)^{\textrm{si}, > 0}$ of $\mathbb{R}\mathbf{Perf}(S)^{\textrm{si}, > 0}$ is a stack whose coarse moduli space $\mathrm{Perf}(S)^{\textrm{si}, > 0}$ is exactly the moduli space Inaba calls $\mathrm{Splcpx}_{S/\mathbb{C}}^{\textrm{\'et}}$ in \cite[\S 3]{ina}.\\

Coming back to Section \ref{projection}, right after Prop. \ref{gammak3}, we consider the projection $\textrm{pr}_{\textrm{der}}$ of $\mathbb{R}\mathbf{Pic}(S)$ onto its full derived factor $$\xymatrix{ \mathbb{R}\mathbf{Pic}(S) \ar[r]^-{\textrm{pr}_{\textrm{der}}} & \mathbb{R}\mathrm{Spec}(\mathrm{Sym} (H^{0}(S, K_{S})[1])).}$$

\begin{df}
The \emph{reduced} derived stack $\mathbb{R}\mathbf{Perf}(S)^{\textrm{si}, \textrm{red}}$ of simple perfect complexes on $S$ is defined by the following homotopy pullback diagram $$\xymatrix{\mathbb{R}\mathbf{Perf}(S)^{\textrm{si}, \textrm{red}} \ar[r] \ar[dd] & \mathbb{R}\mathbf{Perf}(S)^{\textrm{si}, > 0} \ar[d]^-{\textrm{det}_{S}} \\  & \mathbb{R}\mathbf{Pic}(S) \ar[d]^-{\textrm{pr}_{\textrm{der}}} \\ \mathrm{Spec}\, \mathbb{C} \ar[r]_-{x_{0}} & \mathbb{R}\mathrm{Spec}(\mathrm{Sym} (H^{0}(S, K_{S})[1]))  }$$
\end{df}

Since the truncation functor commutes with homotopy pullbacks, the truncation of $\mathbb{R}\mathbf{Perf}(S)^{\textrm{si}, \textrm{red}}$ is
the same as the truncation of $\mathbb{R}\mathbf{Perf}(S)^{\textrm{si}, > 0}$, i.e. $\mathbf{Perf}(S)^{\textrm{si}, > 0}$, therefore its coarse moduli space is again Inaba's $\mathrm{Splcpx}_{S/\mathbb{C}}^{\textrm{\'et}}$ (\cite[\S 3]{ina}).

\begin{thm}\label{inaba}
The composite map $$\xymatrix{\mathbb{R}\mathbf{Perf}(S)^{\textrm{si}, > 0} \ar[r]^-{\textrm{det}_{S}}   & \mathbb{R}\mathbf{Pic}(S) \ar[r]^-{\textrm{pr}_{\textrm{der}}} & \mathbb{R}\mathrm{Spec}(\mathrm{Sym} (H^{0}(S, K_{S})[1])) }$$ is smooth. Therefore the derived stack $\mathbb{R}\mathbf{Perf}(S)^{\textrm{si}, \textrm{red}}$ is actually a smooth, usual (i.e. underived) stack, and $$\mathbb{R}\mathbf{Perf}(S)^{\textrm{si}, \textrm{red}}\simeq \mathrm{t}_{0}(\mathbb{R}\mathbf{Perf}(S)^{\textrm{si}, \textrm{red}}) \simeq \mathbf{Perf}(S)^{\textrm{si}, > 0}.$$ Under these identifications, the canonical map  $\mathbb{R}\mathbf{Perf}(S)^{\textrm{si}, \textrm{red}} \rightarrow \mathbb{R}\mathbf{Perf}(S)^{\textrm{si}, > 0}$ becomes isomorphic to the inclusion of the truncation $\mathbf{Perf}(S)^{\textrm{si}, >0} \rightarrow \mathbb{R}\mathbf{Perf}(S)^{\textrm{si}, > 0}$.
\end{thm}

\noindent \textbf{Proof.} Let $E$ be a perfect complex on $S$ such that $\mathrm{Ext}_{S}^i(E,E)=0$ for $i < 0$, and the trace map $\mathrm{Ext}_{S}^0(E,E)\rightarrow H^{0}(S,\mathcal{O}_{S})\simeq \mathbb{C}$ is an isomorphism. The homotopy fiber product defining $\mathbb{R}\mathbf{Perf}(S)^{\textrm{si}, \textrm{red}}$ yields a distinguished triangle of tangent complexes $$\xymatrix{ \mathbb{T}_{E}\mathbb{R}\mathbf{Perf}(S)^{\textrm{si}, \textrm{red}} \ar[r] & \mathbb{T}_{E}\mathbb{R}\mathbf{Perf}(S)^{\textrm{si}, >0} \ar[r] & H^{0}(S, K_S)^{\vee}[-1]}.$$ Since $$\mathbb{T}_{E}\mathbb{R}\mathbf{Perf}(S)^{\textrm{si}, >0} \simeq \mathbb{R}\mathrm{End}_{S}(E)[1],$$ this complex is cohomologically concentrated in degrees $[-1,1]$. Therefore, to prove the theorem, it is enough to show that the map (induced by the above triangle on $H^{1}$) $$\alpha: \mathrm{Ext}_{S}^2(E,E)\simeq H^{1}(\mathbb{T}_{E}\mathbb{R}\mathbf{Perf}(S)^{\textrm{si}, >0}) \longrightarrow H^{0}(S, K_S)^{\vee}$$ is an isomorphism. If we denote by $$\xymatrix{\alpha' : \mathrm{Ext}_{S}^2(E,E) \ar[r]^-{\alpha} & H^{0}(S, K_S)^{\vee} \ar[r]_-{\sim}^-{\textrm{s}} & H^{2}(S,\mathcal{O}_{S})}$$ (the isomorphism s given by Serre duality),  the following diagram $$\xymatrix{\mathrm{Ext}_{S}^2(E,E) \ar[r]^-{\textrm{s}} \ar[d]_-{\alpha'} & \mathrm{Ext}_{S}^0(E,E)^{\vee} \ar[d]^-{tr_{E}^{\vee}} \\ H^2(S, \mathcal{O}_{S}) \ar[r]_-{\textrm{s}} & H^0(S, \mathcal{O}_{S})^{\vee} }$$
(where again, the s isomorphisms are given by Serre duality on $S$) is commutative. But, by hypothesis, $\textrm{tr}_{E}$ is an isomorphism and we conclude.





\hfill $\Box$ \\

\noindent The following corollary was first proved by Inaba \cite[Thm. 3.2]{ina}.

\begin{cor}
The coarse moduli space $\mathrm{Perf}(S)^{\textrm{si}, > 0}$ of simple perfect complexes on a smooth projective $K3$ surface $S$ is a smooth algebraic space.
\end{cor}

\noindent \textbf{Proof.} The stack  $\mathbb{R}\mathbf{Perf}(S)^{\textrm{si}, \textrm{red}}\simeq \mathbf{Perf}(S)^{\textrm{si}, > 0}$ is a $\mathbb{G}_{m}$-gerbe so its smoothness is equivalent to the smoothness of its coarse moduli space. \hfill $\Box$ \\

\begin{rmk} \emph{Inaba shows in \cite[Thm. 3.3 ]{ina} (again generalizing earlier results by Mukai in \cite{muk}) that the coarse moduli space $\mathrm{Perf}(S)^{\textrm{si}, > 0}$ also carries a canonical \emph{symplectic structure}. In \cite{ptvv} it is shown that in fact the whole derived stack $\mathbb{R}\textbf{Perf}(S)$ carries a natural \emph{derived symplectic structure} of degree $0$, and that this induces on $\mathrm{Perf}(S)^{\textrm{si}, > 0}$ the symplectic structure defined by Inaba.}
\end{rmk}

\begin{rmk} \emph{The same argument used in proving Theorem \ref{inaba} works by replacing $S$ by the de Rham moduli space $\textbf{M}_\textbf{DR}(C)$,  the Dolbeault moduli space $\textbf{M}_\textbf{Dol}(C)$, for $C$ a complex smooth projective curve, or by the Betti moduli space $\textbf{M}_\textbf{B}(S)$ representing the homotopy type of an oriented
compact topological surface $S$ (see \cite{sim}, for the definitions of these moduli spaces). This yields smoothness of $\mathrm{Perf}(X)^{\textrm{si}, > 0}$ for $X=\textbf{M}_\textbf{DR}(C), \textbf{M}_\textbf{Dol}(C), \textbf{M}_\textbf{B}(S)$.}
\end{rmk}

\subsection{On Calabi-Yau threefolds}\label{gwdt}

In this section, for $X$ an arbitrary complex smooth projective variety, we first elaborate on the map $$\mathrm{A}^{(n)}_{X}: \mathbb{R}\overline{\mathbf{M}}_{g,n}(X) \longrightarrow \mathbb{R}\mathbf{Perf}(X)$$ from Def.  \ref{delta1pointed}. This elaboration will give us a  map $\mathrm{C}^{(n)}_{X, \mathcal{L}}$ from a derived substack of $\mathbb{R}\overline{\mathbf{M}}_{g,n}(X)$ to the derived stack $\mathbb{R}\mathbf{Perf}(X)^{\textrm{si}, > 0}_{\mathcal{L}}$, $\mathcal{L}$ being a line bundle on $X$ (see Def. \ref{simples}). When we specialize to the case where $X$ is a projective smooth \emph{Calabi-Yau $3$-fold} $Y$, we prove that $\mathbb{R}\mathbf{Perf}(Y)^{\textrm{si}, > 0}_{\mathcal{L}}$ is \emph{quasi-smooth} (Proposition \ref{adt}), and that the map $\mathrm{C}^{(n)}_{Y, \mathcal{L}}$ allows us to compare the induced obstruction theories on the truncations of its source and target.\\

To begin with, let $X$ be a smooth complex projective variety. 
First of all, observe that taking tensor products of complexes induces an action of the derived group stack $\mathbb{R}\mathbf{Pic}(X)$ on $\mathbb{R}\mathbf{Perf}(X)$ $$\mu:\mathbb{R}\mathbf{Pic}(X) \times \mathbb{R}\mathbf{Perf}(X) \longrightarrow \mathbb{R}\mathbf{Perf}(X).$$ Let  $x_{\mathcal{L}}: \mathrm{Spec}\, \mathbb{C} \rightarrow \mathbb{R}\mathbf{Pic}(X)$ be the point corresponding to a line bundle $\mathcal{L}$ on $X$. 

\begin{df} Let $\sigma_{\mathcal{L}}:  \mathbb{R}\mathbf{Pic}(X) \rightarrow \mathbb{R}\mathbf{Pic}(X)$ be the composite $$\xymatrix{\mathbb{R}\mathbf{Pic}(X) \ar[r]^-{(\mathrm{inv}, x_{\mathcal{L}})} & \mathbb{R}\mathbf{Pic}(X) \times \mathbb{R}\mathbf{Pic}(X)   \ar[r]^-{\times} & \mathbb{R}\mathbf{Pic}(X)} $$ where $\times$ (resp. $\mathrm{inv}$) denotes the product (resp. the inverse) map in $\mathbb{R}\mathbf{Pic}(X)$ (shortly, $\sigma_{\mathcal{L}}(\mathcal{L}_1)= \mathcal{L}\otimes \mathcal{L}^{-1}_{1}$).
\begin{itemize}
\item Define $\mathrm{A}^{(n)}_{X, \mathcal{L}}: \mathbb{R}\overline{\mathbf{M}}_{g,n}(X) \rightarrow \mathbb{R}\mathbf{Perf}(X)_{\mathcal{L}} $ via the composite $$\xymatrix{\mathbb{R}\overline{\mathbf{M}}_{g,n}(X) \ar[rr]^-{(det\, \circ \mathrm{A}^{(n)}_{X}, \mathrm{A}^{(n)}_{X})} & & \mathbb{R}\mathbf{Pic}(X) \times \mathbb{R}\mathbf{Perf}(X) \ar[r]^-{\sigma_{\mathcal{L}} \times \mathrm{id}} & \mathbb{R}\mathbf{Pic}(X) \times \mathbb{R}\mathbf{Perf}(X) \ar[r]^-{\mu} & \mathbb{R}\mathbf{Perf}(X)}$$ (shortly, $\mathrm{A}^{(n)}_{X, \mathcal{L}}(E)= E\otimes (\det{E})^{-1}\otimes \mathcal{L}$).
\item Define the derived open substack $\mathbb{R}\overline{\mathbf{M}}_{g,n}(X)^{emb}\hookrightarrow \mathbb{R}\overline{\mathbf{M}}_{g,n}(X)$ as $\phi_{\mathbb{R}\overline{\mathbf{M}}_{g,n}(X)}(\overline{\mathbf{M}}_{g,n}(X)^{emb})$ (see Prop. \ref{opens}) where $\overline{\mathbf{M}}_{g,n}(X)^{emb}$ is the open substack of $\overline{\mathbf{M}}_{g,n}(X)$ consisting of pointed stable maps which are closed immersions.
\item  Define $\mathrm{C}^{(n)}_{X, \mathcal{L}}: \mathbb{R}\overline{\mathbf{M}}_{g,n}(X)^{emb} \rightarrow \mathbb{R}\mathbf{Perf}(X)^{\textrm{si}, > 0}_{\mathcal{L}}$ via the composite $$\xymatrix{\mathbb{R}\overline{\mathbf{M}}_{g,n}(X)^{emb} \ar@{^{(}->}[r] & \mathbb{R}\overline{\mathbf{M}}_{g,n}(X) \ar[r]^-{\mathrm{A}^{(n)}_{X, \mathcal{L}}} & \mathbb{R}\mathbf{Perf}(X)_{\mathcal{L}}}$$ (note that this composite indeed factors through $\mathbb{R}\mathbf{Perf}(X)^{\textrm{si}, > 0}_{\mathcal{L}}$, since $$\mathrm{tr}:\mathrm{Ext}^{0}(\mathbb{R}f_{*}\mathcal{O}_{C},\mathbb{R}f_{*}\mathcal{O}_{C})\simeq \mathbb{C},$$ if the pointed stable map $f$ is a closed immersion).
\end{itemize}
\end{df}

\begin{rmk} \emph{The map $\mathrm{C}^{(n)}_{X, \mathcal{L}}$ is also defined on the a priori larger open derived substack consisting (in the sense of Prop. \ref{opens}) of pointed stable maps $f$ such that the trace map $\mathrm{tr}:\mathrm{Ext}^{0}(\mathbb{R}f_{*}\mathcal{O}_{C},\mathbb{R}f_{*}\mathcal{O}_{C})\rightarrow H^{0}(X,\mathcal{O}_{X})\simeq \mathbb{C}$ is an isomorphism.} \end{rmk}

We would like to use the map $\mathrm{C}^{(n)}_{X, \mathcal{L}}$ to induce a comparison map between the induced obstruction theories on the truncations of $\mathbb{R}\overline{\mathbf{M}}_{g,n}(X)^{emb}$ and of  $\mathbb{R}\mathbf{Perf}(X)^{\textrm{si}, > 0}_{\mathcal{L}}$.\\
This is possible when we take $X$ to be a Calabi-Yau $3$-fold Y. In fact:\\

\begin{thm}\label{adt} If $Y$ is a smooth complex projective Calabi-Yau $3$-fold, then the derived stack $\mathbb{R}\mathbf{Perf}(Y)^{\textrm{si}, > 0}_{\mathcal{L}}$ is quasi-smooth. Therefore, the closed immersion $j: \mathbf{Perf}(Y)^{\textrm{si}, > 0}_{\mathcal{L}} \hookrightarrow \mathbb{R}\mathbf{Perf}(Y)^{\textrm{si}, > 0}_{\mathcal{L}}$ induces a $[-1,0]$-perfect obstruction theory $j^{*}\mathbb{T}_{\mathbb{R}\mathbf{Perf}(Y)^{\textrm{si}, > 0}_{\mathcal{L}}} \rightarrow \mathbb{T}_{\mathbf{Perf}(Y)^{\textrm{si}, > 0}_{\mathcal{L}}}$.
\end{thm}

\noindent \textbf{Proof.} This is a corollary of Proposition \ref{almostDT}. Let $\mathbb{T}_{E}$ be the tangent complex of $\mathbb{R}\mathbf{Perf}(Y)^{\textrm{si}, > 0}_{\mathcal{L}}$ at a point corresponding to the perfect complex $E$. Now $Y$ is Calabi-Yau of dimension $3$, so $\Omega_{Y}^3 \equiv K_{Y} \simeq \mathcal{O}_Y$; but $E$ is simple (i.e. the trace map $\mathrm{Ext}^0(F,F)\rightarrow H^{0}(X,\mathcal{O}_{X})\simeq \mathbb{C}$ is an isomorphism), so Serre duality implies $\mathrm{Ext}^i(E,E)_0 =0$ for $i\geq 3$ (and all $i\leq 0$). Therefore the perfect complex $\mathbb{T}_{E}$ is concentrated in degrees $[0,1]$, and $\mathbb{R}\mathbf{Perf}(Y)^{\textrm{si}, > 0}_{\mathcal{L}}$ is quasi-smooth. The second assertion follows immediately from Prop. \ref{inducedobstruction}.
\hfill $\Box$ \\

\begin{rmk} \emph{Note that the stack $\mathbf{Perf}(Y)^{\textrm{si}, > 0}_{\mathcal{L}}$ is \emph{not} proper over $\mathrm{Spec} \, \mathbb{C}$. However it receives maps from both Thomas moduli space $I_{n}(Y;\beta)$ of ideal sheaves (whose subschemes have Euler characteristic $n$ and fundamental class $\beta \in H_{2}(Y,\mathbb{Z})$) - see \cite{tho} -  and from Pandharipande-Thomas moduli space $P_{n}(Y; \beta)$ of stable pairs - see \cite{pt}. For example, the map from $P_{n}(Y;\beta)$ sends a pair to the pair itself, considered as a complex on $Y$. Moreover, at the points in the image of such maps, the tangent and obstruction spaces of these spaces are the same as those induced from the cotangent complex of our $\mathbb{R}\mathbf{Perf}(Y)^{\textrm{si}, > 0}_{\mathcal{L}}$ (see e.g. \cite[\S 2.1]{pt}).} 
\end{rmk}
 
As showed in \S \ref{functobs}, the map $\mathrm{C}^{(n)}_{Y, \mathcal{L}}: \mathbb{R}\overline{\mathbf{M}}_{g,n}(Y)^{emb} \longrightarrow \mathbb{R}\mathbf{Perf}(Y)^{\textrm{si}, > 0}_{\mathcal{L}}$ induces a comparison map between the two obstruction theories. More precisely, the commutative diagram in $\dst$ $$\xymatrix{\overline{\mathbf{M}}_{g,n}(Y)^{emb} \ar[r]^-{\tr \mathrm{C}^{(n)}_{Y, \mathcal{L}}} \ar@{^{(}->}[d]_-{j_{GW}} &  \mathbf{Perf}(Y)^{\textrm{si}, > 0}_{\mathcal{L}} \ar@{^{(}->}[d]^-{j_{DT}} \\  \mathbb{R}\overline{\mathbf{M}}_{g,n}(Y)^{emb} \ar[r]_-{\mathrm{C}^{(n)}_{Y, \mathcal{L}}} & \mathbb{R}\mathbf{Perf}(Y)^{\textrm{si}, > 0}_{\mathcal{L}} }$$(where each $j$ is the closed immersion of the truncation of a derived stack into the full derived stack), induces a morphism of triangles $$\xymatrix{(\tr \mathrm{C}^{(n)}_{Y, \mathcal{L}})^*j_{DT}^*\mathbb{L}_{\mathbb{R}\mathbf{Perf}(Y)^{\textrm{si}, > 0}_{\mathcal{L}}} \ar[r] \ar[d] & (\tr \mathrm{C}^{(n)}_{Y, \mathcal{L}})^*\mathbb{L}_{\mathbf{Perf}(Y)^{\textrm{si}, > 0}_{\mathcal{L}}} \ar[r] \ar[d] & (\tr \mathrm{C}^{(n)}_{Y, \mathcal{L}})^*\mathbb{L}_{\mathbb{R}\mathbf{Perf}(Y)^{\textrm{si}, > 0}_{\mathcal{L}}/\mathbf{Perf}(Y)^{\textrm{si}, > 0}_{\mathcal{L}}} \ar[d] \\ j_{GW}^*\mathbb{L}_{\mathbb{R}\overline{\mathbf{M}}_{g,n}(Y)^{emb}} \ar[r] & \mathbb{L}_{\overline{\mathbf{M}}_{g,n}(Y)^{emb}} \ar[r] & \mathbb{L}_{\mathbb{R}\overline{\mathbf{M}}_{g,n}(Y)^{emb}/\overline{\mathbf{M}}_{g,n}(Y)^{emb}} }$$
- in the derived category of perfect complexes on $\overline{\mathbf{M}}_{g,n}(Y)^{emb}$ - i.e. a morphism relating the two obstruction theories induced on the truncations stacks $\overline{\mathbf{M}}_{g,n}(Y)^{emb}$ and $\mathbf{Perf}(X)^{\textrm{si}, > 0}_{\mathcal{L}}$. Note that, for the object in the upper left corner of the above diagram, we have a natural isomorphism $$(\tr \mathrm{C}^{(n)}_{Y, \mathcal{L}})^*j_{DT}^*\mathbb{L}_{\mathbb{R}\mathbf{Perf}(Y)^{\textrm{si}, > 0}_{\mathcal{L}}} \simeq j_{GW}^*(\mathrm{C}^{(n)}_{Y, \mathcal{L}})^*\mathbb{L}_{\mathbb{R}\mathbf{Perf}(Y)^{\textrm{si}, > 0}_{\mathcal{L}}}.$$ \\

\begin{rmk} \emph{The motivation for constructing the above comparison morphism between the induced obstruction theories comes from the so called \emph{GW/DT comparison}. The conjectural comparison between Gromov-Witten and Donaldson-Thomas or Pandharipande-Thomas  invariants for a general Calabi-Yau $3$-fold $Y$ (stated in \cite{mnop}, and proved in special cases, e.g. \cite{bb} and \cite{moop}) might be in principle approached by producing a map from the stack of pointed stable maps to $Y$ to Thomas moduli space of ideal sheaves (\cite{tho}) or to the Pandharipande-Thomas stack of stable pairs on $Y$ (\cite{pt}), \emph{together} with an induced map relating the corresponding obstruction theories. It is worth remarking that simply producing a map between these moduli spaces is \emph{not} enough to relate the standard obstruction theories (the ones yielding the GW and DT or PT invariants, respectively): the notion of obstruction theory is not canonically attached to a stack, and therefore has not enough functoriality. This problem is completely solved in when the obstruction theories come from derived extensions, as explained in \S \ref{derobs}. Therefore, one could hope to produce instead a map from the derived stack of pointed stable maps to $Y$ to a derived extension of the moduli space of ideal sheaves on $Y$ or of the stack of stable pairs on $Y$. We are not able to do this, at present. What we did above was to produce a map of derived stacks from  the derived open substack $\mathbb{R}\overline{\mathbf{M}}_{g,n}(Y)^{emb}$ of $\mathbb{R}\overline{\mathbf{M}}_{g,n}(Y)$ consisting of stable maps which are closed immersions, to the derived stack $\mathbb{R}\mathbf{Perf}(Y)^{\textrm{si}, > 0}_{\mathcal{L}}$ of simple perfect complexes with no negative Ext's and fixed determinant $\mathcal{L}$ (for any $\mathcal{L}$). As shown above, such a map automatically induces a morphism between the corresponding obstruction theories, with no need of further data, and one might think of this as a baby, open version of the GW/DT comparison.\\ In the presence of suitable compactifications $\mathbb{R}\mathcal{N}$ of $\mathbb{R}\overline{\mathbf{M}}_{g,n}(Y)^{emb}$ and  $\mathbb{R}\mathcal{P}$ of $\mathbb{R}\mathbf{Perf}(Y)^{\textrm{si}, > 0}_{\mathcal{L}}$ (or better of some derived substacks thereof, cut out by suitable cohomological or $K$-theoretic numerical conditions), to which the map $\mathrm{C}^{(n)}_{Y, \mathcal{L}}$ extends as a quasi-smooth map $F: \mathbb{R}\mathcal{N} \rightarrow \mathbb{R}\mathcal{P}$, the corresponding morphism between obstruction theories would give a canonical way of comparing the induced virtual fundamental classes, and therefore, the two counting invariants. More precisely, if $f$ denotes the truncation of the morphism $F$, we would get (\cite[Thm. 7.4]{sch} or Prop. \ref{timo}) $$f^!([\mathcal{P}]^\textrm{vir})=[\mathcal{N}]^\textrm{vir},$$ where $f^!:A_*(\mathcal{P}) \rightarrow A_*(\mathcal{N})$ denotes the \emph{virtual pullback} between Chow groups, as defined in \cite{manolache}.\\
Of course the real technical heart of the GW/DT comparison lies exactly in the fine analysis of what happens at the boundary of the compactifications of the two stacks involved (in the above picture, the existence of an extension $F: \mathbb{R}\mathcal{N} \rightarrow \mathbb{R}\mathcal{P}$), and for this our methods from derived algebraic geometry does not supply, at the moment, any new tool or direction.}
\end{rmk}

\begin{rmk} \emph{The group $\mathrm{Aut}_{\textrm{Perf}}(X)$ of self-equivalences of the derived category $\mathrm{D}_{\textrm{Perf}}(X)$ of perfect complexes on $X$, acts on the derived stack $\mathbb{R}\mathbf{Perf}(X)$. This action preserves, in an obvious sense, the induced obstruction theory on the stack $\mathbf{Perf}(X)$. It would be interesting to study how this action affects Bridgeland stability conditions on $\mathrm{D}_{\textrm{Perf}}(X)$, and use it to draw consequences on the comparison of various counting invariants (e.g. Donaldson-Thomas and Pandharipande-Thomas ones). This will be the object of a future paper.}
\end{rmk}

\begin{appendix}

\section{Derived stack of perfect complexes and Atiyah classes}\label{AXexplained}

We explain here the relationship between the tangent maps associated to morphisms to the stack of perfect complexes and Atiyah classes (of perfect complexes) used in the main text (see \S \ref{AX}). As in the main text, we work over $\mathbb{C}$, even if most of what we say below holds true over any field of characteristic zero. As usual, all tensor products and fiber products will be implicitly derived.\\

If $\mathcal{Y}$ is a derived geometric stack having a perfect cotangent complex (\cite[\S 1.4]{hagII}), and $E$ is a perfect complex on $\mathcal{Y}$, then we will implicitly identify the Atiyah class map of $E$  $$\mathrm{at}_{E}: E \longrightarrow \mathbb{L}_{\mathcal{Y}} \otimes E[1] $$  with the corresponding map $$\mathbb{T}_{\mathcal{Y}} \longrightarrow E^{\vee}\otimes E[1]$$ via the bijection $$[\mathbb{T}_{\mathcal{Y}},E^{\vee}\otimes E [1]]\simeq [\mathbb{T}_{\mathcal{Y}}\otimes E, E[1]]\simeq [E,\mathbb{L}_{\mathcal{Y}}\otimes E[1]]$$
given by the adjunction $(\otimes, \mathbb{R}\underline{\mathrm{Hom}})$, and perfectness of $E$ and $\mathbb{L}_{\mathcal{Y}}$ (where $[-,-]$ denotes the $Hom$ set in the derived category of perfect complexes on $\mathcal{Y}$).\\

We start with a quite general situation. Let $\mathcal{Y}$ be a derived geometric stack having a perfect cotangent complex, and $\mathbf{Perf}$ the stack of perfect complexes (viewed as a derived stack). Then, giving a map of derived stacks $\phi_{E}:\mathcal{Y}\rightarrow \mathbf{Perf}$ is the same thing as giving a perfect complex $E$ on $\mathcal{Y}$, and 
\begin{itemize}
\item $\phi_{E}^{*} \mathbb{T}_{\mathbf{Perf}} \simeq \mathbb{R}\underline{\mathrm{End}}_{\mathcal{Y}}(E)[1]$
\item the tangent map to $\phi_{E}$ $$\xymatrix{\mathbb{T}\phi_{E}:\mathbb{T}_{\mathcal{Y}} \ar[rr] & &  \phi_{E}^{*} \mathbb{T}_{\mathbf{Perf}} \simeq \mathbb{R}\underline{\mathrm{End}}_{\mathcal{Y}}(E)[1]\simeq E^{\vee}\otimes E [1]}$$ is the Atiyah class map $\mathrm{at}_{E}$ of $E$.
\end{itemize}

\begin{rmk}\emph{The second point above might be considered as a definition when $\mathcal{Y}$ is a derived stack, and it coincides with Illusie's definition (\cite[Ch. 4, 2.3.7]{ill}) when $\mathcal{Y}= Y$ is a quasi-projective scheme. In fact, in this case, the map $\Phi_{E}$ factors through the stack $\mathbf{Perf}^{\textrm{strict}}$ of strict perfect complexes; thus the proof reduces immediately to the case where $E$ is a vector bundle on $Y$, which is straightforward.} 
\end{rmk}

\noindent The above description applies in particular to a map of derived stacks of the form $$\Phi_{E}:\mathcal{Y}:=\mathcal{X} \times X \longrightarrow \mathbf{Perf}$$ where $X$ is a smooth projective scheme, $\mathcal{X}$ is a derived geometric stack having a perfect cotangent complex, and $E$ is a perfect complex on $\mathcal{X} \times X$ : in the main text we are interested in $\mathcal{X}=\mathbb{R}\overline{\mathbf{M}}_{g}(X)$.
Such a map corresponds, by adjunction to a map $$\Psi_{E}:\mathcal{X}\longrightarrow \mathbb{R}\mathrm{HOM}(X, \mathbf{Perf})=: \mathbb{R}\mathbf{Perf}(X).$$ The tangent map of $\Psi_{E}$ fits into the following commutative diagram $$\xymatrix{\mathbb{T}_{\mathcal{X}} \ar[r]^-{\mathbb{T}\Psi_{E}} \ar[d]_-{\textrm{can}} & \Psi_{E}^{*}\mathbb{T}_{\mathbb{R}\mathbf{Perf}(X)} \ar[r]^-{\sim} & \mathbb{R}pr_{\mathcal{X},*}(E^{\vee}\otimes E)[1] \\ \mathbb{R}pr_{\mathcal{X},*}pr_{\mathcal{X}}^{*}\mathbb{T}_{\mathcal{X}} \ar[r]_-{\textrm{can}} & \mathbb{R}pr_{\mathcal{X},*}(pr_{\mathcal{X}}^{*}\mathbb{T}_{\mathbb{X}}\oplus pr_{X}^{*}\mathbb{T}_{X})  \ar[r]_-{\sim} & \mathbb{R}pr_{\mathcal{X},*}\mathbb{T}_{\mathcal{X} \times X} \ar[u]_-{\mathbb{R}pr_{\mathcal{X},*}(\mathbb{T}\Phi_{E})} }
$$ where $can$ denote obvious canonical maps, and we can identify $\mathbb{R}pr_{\mathcal{X},*}(\mathbb{T}\Phi_{E})$ with $\mathbb{R}pr_{\mathcal{X},*}(\mathrm{at}_{E})$, in the sense explained above. In other words, $\mathbb{T}\Psi_{E}$ is described in terms of the relative Atiyah class map $$\mathrm{at}_{E/X}: pr_{\mathcal{X}}^{*}\mathbb{T}_{\mathcal{X}}\simeq \mathbb{T}_{\mathcal{X} \times X/X}\rightarrow E^{\vee}\otimes^{\mathbb{L}} E [1]$$ of $E$ relative to $X$, as the composition $$\xymatrix{\mathbb{T}\Psi_{E}:\mathbb{T}_{\mathcal{X}}\ar[r]^-{\textrm{can}} & \mathbb{R}pr_{\mathcal{X},*}pr_{\mathcal{X}}^{*}\mathbb{T}_{\mathcal{X}} \ar[r]^-{\sim} &  \mathbb{R}pr_{\mathcal{X},*}\mathbb{T}_{\mathcal{X} \times X /X} \ar[rrr]^-{\mathbb{R}pr_{\mathcal{X},*}(\mathrm{at}_{E/X})} & & &  \mathbb{R}pr_{\mathcal{X},*}(E^{\vee}\otimes E)[1] . }$$ 
\begin{rmk}\emph{$\mathbb{T}\Psi_{E}$ might be viewed at as a generalization of what is sometimes called the Kodaira-Spencer map associated to the $\mathcal{X}$-family $E$ of perfect complexes over $X$ (e.g. \cite[formula (14)]{km}).}
\end{rmk}

In the main text, we are interested in the case $\mathcal{X}=\mathbb{R}\overline{\mathbf{M}}_{g}(X)$, $\mathrm{pr}:=\mathrm{pr}_{\mathcal{X}}$, and $E$ perfect of the form $\mathbb{R}\pi_{*}\mathcal{E}$, where $$\pi: \mathbb{R}\mathcal{C}_{g;\,X}\longrightarrow \mathbb{R}\overline{\mathbf{M}}_{g}(X)\times X$$ is the universal map and $\mathcal{E}$ is a complex on $\mathbb{R}\mathcal{C}_{g;\,X}$, namely $\mathcal{E}=\mathcal{O}_{\mathbb{R}\mathcal{C}_{g;\,X}}$. In such cases, if we call $(f:C\rightarrow X)$ the stable map corresponding to the complex point $x$, we have a ladder of homotopy cartesian diagrams 
$$\xymatrix{C \ar[r]^-{\iota_{f}} \ar[d]_{f} & \mathbb{R}\mathcal{C}_{g;\,X} \ar[d]^-{\pi}  & \\ X \ar[r]^-{\underline{x}} \ar[d]_-{q} & \mathbb{R}\overline{\mathbf{M}}_{g}(X)\times X \ar[d]^-{\mathrm{pr}} \ar[r]^-{\mathrm{pr}_{X}} & X \ar[d]^-{q} \\ \mathrm{Spec}\, \mathbb{C} \ar[r]_-{x} & \mathbb{R}\overline{\mathbf{M}}_{g}(X) \ar[r] & \mathrm{Spec}\, \mathbb{C}}$$ and the base-change isomorphism (true in derived algebraic geometry with no need of flatness) gives us $$\underline{x}^{*}E=\underline{x}^{*}\mathbb{R}\pi_{*}\mathcal{E} \simeq \mathbb{R}f_{*}\iota_{f}^{*}\mathcal{E}.$$ For $\mathcal{E}=\mathcal{O}_{\mathbb{R}\mathcal{C}_{g;\,X}}$, we then get $$\underline{x}^{*}E=\underline{x}^{*}\mathbb{R}\pi_{*}\mathcal{O}_{\mathbb{R}\mathcal{C}_{g;\,X}} \simeq \mathbb{R}f_{*}\mathcal{O}_{C}.$$ Again by base-change formula, we get $x^*\mathbb{R}\mathrm{pr}_{*}\simeq \mathbb{R}q_{*}\underline{x}^{*}$, and therefore the tangent map to $\mathrm{A}_{X}:=\Psi_{\mathbb{R}\pi_{*}\mathcal{O}_{\mathbb{R}\mathcal{C}_{g;\,X}}}$ at $x=(f:C\rightarrow X)$, is the composition $$\xymatrix{\mathbb{T}_{x}\mathrm{A}_{X}: \mathbb{T}_{x} \mathbb{R}\overline{\mathbf{M}}_{g}(X) \simeq \mathbb{R}\Gamma(C, \mathrm{Cone}(\mathbb{T}_{C}\rightarrow f^{*}T_{X})) \ar[r] & \mathbb{R}\Gamma(X, \underline{x}^{*}\mathbb{T}_{\mathbb{R}\overline{\mathbf{M}}_{g}(X) \times X}) \ar[r]  & }$$ $$\xymatrix{ \ar[rrr]^-{\mathbb{R}\Gamma(X,\underline{x}^{*}\textrm{at}_{E})} & & & \mathbb{R}\mathrm{End}_{X}(\mathbb{R}f_{*}\mathcal{O}_{C})[1] \simeq \mathbb{T}_{\mathbb{R}f_*\mathcal{O}_{C}}\mathbb{R}\mathbf{Perf}(X)}$$ 
\medskip

The following is the third assert in Proposition \ref{finallyatiyah}, \S \ref{AX}.

\begin{prop} \emph{The composition $$\xymatrix{\mathbb{R}\Gamma(X,T_X) \ar[r]^-{\textrm{can}} & \mathbb{R}\Gamma(X,\mathbb{R}f_{*}f^{*}T_X) \ar[r]^-{\textrm{can}} &   \mathbb{R}\Gamma(X,\textrm{Cone}(\mathbb{R}f_{*}\mathbb{T}_{C}\rightarrow \mathbb{R}f_{*}f^{*}T_X))\simeq \mathbb{T}_{x} \mathbb{R}\overline{\mathbf{M}}_{g}(X) \ar[r] & }$$ $$\xymatrix{ \ar[r]^-{\mathbb{T}_{x}\mathrm{A}_{X}} & x^{*}\mathrm{A}_{X}^{*}\mathbb{T}\mathbb{R}\mathbf{Perf}(X)\simeq \mathbb{T}_{\mathbb{R}f_{*}\mathcal{O}_{C}}\mathbb{R}\mathbf{Perf}(X) \simeq \mathbb{R}\mathrm{End}_{X}(\mathbb{R}f_{*}\mathcal{O}_{C})[1]  }$$ coincides with $\mathbb{R}\Gamma(X, \textrm{at}_{\mathbb{R}f_{*}\mathcal{O}_{C}})$.}
\end{prop}

\noindent\textbf{Proof.} We first observe that if $\mathcal{F}$ is perfect complex on $X$, and $\mathbb{R}\mathbf{Aut}(X)$ is the derived stack of automorphisms of $X$, there are obvious maps of derived stacks $$\rho_{x}: \mathbb{R}\mathbf{Aut}(X)\longrightarrow \mathbb{R}\mathbf{HOM}_{\dst}(C,X)$$  $$\sigma_{\mathcal{F}}:\mathbb{R}\mathbf{Aut}(X)\longrightarrow \mathbb{R}\mathbf{Perf}(X)$$ induced by the natural action of $\mathbb{R}\mathbf{Aut}(X)$ by composition on maps and by pullbacks on perfect complexes, respectively. Moreover, the tangent map to $\sigma_{\mathcal{F}}$ at the identity $\mathrm{Spec}\, \mathbb{C}$-point of $\mathbb{R}\mathbf{Aut}(X)$ $$\mathbb{T}_{\mathrm{id_{X}}}\sigma_{\mathcal{F}}:\mathbb{R}\Gamma(X,T_{X})\simeq \mathbb{T}_{\mathrm{id_{X}}}\mathbb{R}\mathbf{Aut}(X) \longrightarrow \mathbb{T}_{\mathcal{F}}\mathbb{R}\mathbf{Perf}(X)\simeq \mathbb{R}\mathrm{End}_{X}(\mathcal{F})[1]$$ is $\mathbb{R}\Gamma(X, \mathrm{at}_{\mathcal{F}})$, where $\mathrm{at}_{\mathcal{F}}$ is the Atiyah class map of $\mathcal{F}$. Then we observe that, by taking $\mathcal{F}:= \underline{x}^{*}\mathbb{R}\pi_{*}\mathcal{O}_{\mathbb{R}\mathcal{C}_{g;\,X}}$ - which is, by base-change formula, isomorphic to $\mathbb{R}f_{*}\mathcal{O}_{C}$ - we get that the composition $$\xymatrix{k_{x}: \mathbb{R}\mathbf{Aut}(X) \ar[r]^-{\rho_{x}} &  \mathbb{R}\mathbf{HOM}_{\dst}(C,X) \ar[r]^-{\textrm{can}} & \mathbb{R}\overline{\mathbf{M}}_{g}(X) \ar[r]^-{\mathrm{A}_{X}} &  \mathbb{R}\mathbf{Perf}(X) }$$ coincides with $\sigma_{\mathcal{F}}$. But the map in the statement of the proposition is just $\mathbb{T}_{\mathrm{id}_{X}}k_{x}$, and we conclude.

\hfill $\Box$ \\

\end{appendix}

\end{document}